\documentclass[11pt,leqno]{article}
\usepackage[english,activeacute]{babel}
\usepackage{color}
\usepackage{amsmath,amsfonts,amssymb,mathrsfs,amsthm}
\newcommand\R{\mathbb{R}}

\newcommand\N{\mathbb{N}}
\newcommand\Z{\mathbb{Z}}
\newcommand\T{\mathbb{T}}
\newcommand\bna{\begin{eqnarray*}}
\newcommand\ena{\end{eqnarray*}}

\newcommand\bnan{\begin{eqnarray}}
\newcommand\enan{\end{eqnarray}}

\newcommand\Xsbt{X_{b,s}^{T}}
\newcommand\Xsb{X_{b,s}}
\newcommand\Zsbt{Z_{b,s}^{T}}

\newcommand{\pt}{\partial _t}
\newcommand{\px}{\partial _x}
\newcommand\nor[2]{\left\|#1\right\|_{#2}}

\newcommand\kdv[1]{\pt #1+#1\px #1+ \px ^3 #1}

\newcommand\Tu{\mathbb{T}}

\newtheorem{thrm}{Theorem}[section]
\newtheorem{rmrk}[thrm]{Remark}
\newtheorem{lmm}[thrm]{Lemma}
\newtheorem{crllr}[thrm]{Corollary}
\newtheorem{prpstn}[thrm]{Proposition}

\begin{document}
\title{\bf Control and Stabilization of the Korteweg-de Vries
Equation on a Periodic Domain}

\date{}

\author{Camille Laurent \\    \small{Laboratoire de  Math\'ematiques} \\
\small{Universit\'e Paris-Sud}\\
\small{B\^{a}timent 425, F-91405 Orsay Cedex, France} \\
\small{email: camille.laurent@math.u-psud.fr}\\
\quad \bigskip \\
Lionel Rosier\\   \small{Institut \'Elie Cartan}\\ \small{UMR 7502
UHP/CNRS/INRIA, B.P. 239}\\ \small{ F-54506 Vand\oe
uvre-l\`es-Nancy Cedex, France}  \\
\small{email: rosier@iecn.u-nancy.fr }\\ \bigskip \small{and}    \medskip\\
Bing-Yu Zhang \\
\small{Department of Mathematical Sciences}\\ \small{University of Cincinnati}\\
\small{Cincinnati, Ohio 45221, USA}\\
\small{email: bzhang@math.uc.edu }  }

\maketitle

\begin{abstract}
This paper aims at completing an earlier work of  Russell and Zhang
\cite{rz-1} to study internal control problems for the distributed
parameter system described by the Korteweg-de Vries equation on a
periodic domain $\T$.  In \cite{rz-1}, Russell and Zhang  showed
that the system is \emph{locally} exactly controllable and
\emph{locally} exponentially stabilizable when the control acts on an
arbitrary nonempty subdomain of $\T$.
In this paper, we show  that the system is in fact {\em globally}
exactly controllable and {\em globally} exponentially stabilizable.
The global exponential stabilizability corresponding to a natural
feedback law is first established with the
aid of certain properties of propagation of compactness and
propagation of regularity in Bourgain spaces for solutions of the
associated linear system.  Then, using a different feedback
law, the resulting closed-loop system is shown to be  locally
exponentially stable with an {\em arbitrarily large} decay rate. A
{\em time-varying} feedback law is further designed to ensure a
global exponential stability with an arbitrary large decay rate.
\end{abstract}

\section{Introduction}

The well-known  Korteweg-de Vries (KdV) equation can be written as
\begin{equation}
\label{kdv1}
\partial _t u + u\partial_x u + \partial _x ^3 u=0,
\end{equation}
where $u=u(x,t)$ denotes a real-valued function of two real
variables $x$ and $t$. The equation was first derived by Korteweg
and de Vries \cite{kdv} in 1895 (or by Boussinesq \cite{bouss} in 1876
\footnote{The interested readers are refereed to a nice article of
de Jager \cite{jager} for the origin of the KdV equation.})
as a model
for propagation of some surface water waves along a channel.
The KdV equation has been intensively studied from
various aspects of both mathematics and physics since the 1960s when
solitons were discovered through solving the KdV equation, and the
inverse scattering method, a so-called nonlinear Fourier transform,
was invented to seek solitons \cite{soliton, miura}. It turns out
that the equation is not only a good model for some water waves but
also a very useful approximation model in nonlinear studies whenever
one wishes to include and balance a weak nonlinearity and weak
dispersive effects \cite{miura}. In particular, the equation is now
commonly accepted as a mathematical model for the unidirectional
propagation of small-amplitude long waves in nonlinear dispersive
systems.

In this paper, we consider the KdV equation posed on the periodic
domain $\T$:
\begin{equation}
\kdv{u} =0, \quad x\in \T, \ t\in \R . \label{1.0}
\end{equation}
 The equation  is known to possess an infinite set of conserved
integral quantities of which  the first two are
 \[ I_1 (t) = \int
_{\T} u(x,t) dx  \] and
\[ I_2 (t) = \int _{\T} u^2 (x,t) dx . \]

\, From the historical origins \cite{kdv, bouss, miura} of the KdV
equation, involving the behavior of water waves in a shallow
channel, it is natural to think of $I_1$  and $I_2$ as expressing
conservation of volume (or mass) and energy, respectively.  The
Cauchy problem for the equation (\ref{1.0}) has been intensively
studied for many years (see  \cite{st, kato, bourgain-1bis, kpv-2}
and the references therein).  The best known result so far
\cite{KappelerTopalov} is that the Cauchy problem is well-posed
in the space $H^s (\T)$ for any $s\geq -1$:

\medskip
\emph{Let $s\geq -1$ and $T>0$ be given. For any $u_0\in H^s (\T)$,
the equation (\ref{1.0}) admits a unique solution $u\in C([0,T]; H^s
(\T))$ satisfying
\[ u(x,0) = u_0 (x) .\]
Moreover, the corresponding solution map $(u_0 \to u$) is continuous
from the space $H^s (\T)$ to the space $C([0,T]; H^s (\T))$.
\footnote{If $s>-\frac12$, this solution map is, in fact, analytic.}
}

\medskip
In this paper we will study the equation (\ref{1.0}) from a control
point of view with a forcing term $f=f(x,t)$ added to the
equation as a control input:
\begin{equation}
\kdv{u} =f, \quad x\in \T, \ t\in \R, \label{1.0-1}
\end{equation}
where $f$ is assumed to be supported in  a given open set $\omega
\subset \T$.
The following exact control problem
and stabilization problem are fundamental in control theory.

\medskip
\noindent {\bf Exact control problem:} \emph{Given an initial state
$u_0$ and a terminal state $u_1$ in a certain space, can one find an
appropriate control input $f$ so that the equation (\ref{1.0-1})
admits a solution $u$ which satisfies $u(.,0)=u_0$ and $u(.,T)=u_1$?}

\medskip
\noindent {\bf Stabilization problem:} \emph{Can one find a feedback
control law: $f=Ku$ so that the resulting closed-loop system
\[ \kdv{u}=Ku,  \quad x\in \T, \ t\in \R ^+\]
is asymptotically stable as $t \to +\infty $?}

\medskip
The problems were first studied by Russell and Zhang for the KdV
equation \cite{rz, rz-1}. In their work, in order to keep the
\emph{mass} $I_1(t)$ conserved, the control input $f(x,t)$ is
chosen to be of the form
\begin{equation}
 f(x,t) = [Gh](x,t) := g(x) \left (h(x,t)-  \int
_{\T}g(y) h(y,t) dy \right ) \label{1.2} \end{equation} where $h$ is
considered as a new control input,  and $g(x)$ is a  given
nonnegative smooth function such that $\{ g > 0 \} = \omega$ and
\[ 2\pi [g] = \int _{\T} g(x) dx =1 .\]
For the chosen $g$, it is easy to see that

\[ \frac{d}{dt} \int _{\T} u(x,t) dx = \int _{\T} f(x,t) dx=0
\qquad \text{for  any} \ t\in \R\] for any solution  $u=u(x,t)$   of
the system
\begin{equation}\label{1.1}
\kdv{u}=Gh ;\end{equation}
thus the \emph{mass} of the system is indeed
conserved.

The following results are due to Russell and Zhang \cite{rz-1}.

\bigskip
\noindent {\bf Theorem A:}  {\em Let $s\geq 0$ and $T>0$ be given.
There exists a $\delta >0$ such that for any $u_0 ,u_1\in H^s (\T)$
with $[u_0]=[u_1]$ satisfying
\[ \| u_0\|_s \leq \delta, \qquad \| u_1\|_s \leq \delta ,\]
one can find a control input $h\in L^2 (0,T; H^s (\T))$ such that
the system (\ref{1.1}) admits a solution $u\in C([0,T]; H^s (\T))$
satisfying
\[ u(x,0) = u_0 (x), \qquad u(x,T)= u_ 1 (x).\]}

\medskip
In order to stabilize the system (\ref{1.1}), Russell and
Zhang employed a simple feedback control law
\begin{equation}
h(x,t) = - G^*u (x,t). \label{1.3}
\end{equation}
The resulting closed-loop system
\begin{equation}
\kdv{u}=-  GG^*u, \quad x\in \T,\  t\in \R. \label{1.4}
\end{equation}
is \emph{locally}  exponentially stable.

\bigskip
\noindent {\bf Theorem B:} {\em Let   $s = 0$ or $s \geq 1$ be
given. There exist positive constants $M, \delta $ and $\gamma$
 such that if $u_0\in H^s (\T)$ satisfies
 \begin{equation}
 \| u_0-[u_0] \| _s \leq \delta
 ,\label{1.5}
 \end{equation}
 then the corresponding solution u of (\ref{1.4})
satisfies \[ \| u(\cdot, t)-[u_0 ]\| _s \leq Me^{-\gamma t} \|
u_0-[u_0] \| _s\] for any $t\geq 0$.  }

\bigskip

Thus one can always find an appropriate control input $h$ to guide
the system (\ref{1.1}) from a given initial state $u_0$ to a
terminal state $u_1$ so long as \emph{their amplitudes are small} and
$[u_0]=[u_1]$. A question arises naturally.

\medskip
\noindent {\bf Question 1:} \emph{Can one still guide the system by
choosing appropriate control input $h$ from a given initial state
$u_0$ to a given terminal state $u_1$ when $u_0$ or $u_1$ have
\emph{large} amplitude?}

\medskip
As for the closed-loop system (\ref{1.4}), its  \emph{small
amplitude} solutions  decay at a uniform exponential rate to the
corresponding constant state $[u_0]$ with respect to the norm in the
space   $H^s (\T)$ as $t\to \infty $. One may ask naturally:

\medskip
\noindent {\bf Question 2:} {\em Do the large amplitude solutions of
the closed-loop system (\ref{1.4}) decay exponentially as $t\to
\infty $?}

\medskip
A further question is:

\medskip
\noindent {\bf Question 3:} {\em For any given number $\lambda
>0$, can we design a linear feedback control law such that the
exponential decay rate of the resulting closed-loop system is
$\lambda $?}

\medskip
One of the main results in this paper is  a positive answer to
Question 1 as given below.
\begin{thrm}\label{main-4}
Let $s\geq 0$,  $R>0$,  and $\mu\in\R$ be given. There exists a time $T>0$
such that if $u_0, u_1 \in H^s (\T)$ with $[u_0]=[u_1]=\mu$ are such that
\[ \| u_0 \|_s \leq R , \qquad \|u_1\|_s \leq R, \]
then one can find a control input $h\in L^2 (0,T; H^s (\T))$ such
that the system (\ref{1.1}) admits a solution $u\in C([0,T]; H^s
(\T))$ satisfying
\[ u(x,0) = u_0 (x), \qquad u(x,T)= u_ 1 (x).\]
\end{thrm}
So the system (\ref{1.1}) is \emph{globally} exactly controllable.

\medskip
As for Question 2, we have the following affirmative answer.
\begin{thrm}
\label{main-1}
Let $s\geq 0 $ and $\mu\in \R$ be given.  There exists a  constant $\kappa
>0$  such that for any $u_0\in H^s (\T)$ with $[u_0]=\mu$,
 the corresponding solution $u$ of the system (\ref{1.4}) satisfies \[
\| u(\cdot, t)-[u_0 ]\| _s \leq \alpha _{s,\mu} (\| u_0 - [u_0] \|_0)e^{-\kappa t}
\| u_0-[u_0] \| _s \qquad \text{ for all  } t\ge 0 ,\]
where $\alpha _{s,\mu} : \R^+ \to \R^+$ is a nondecreasing  continuous function
depending on $s$ and $\mu$.
\end{thrm}

Note that Theorem \ref{main-4} follows from Theorem \ref{main-1}
and a local control result around the state $u(x)=\mu$ (similar to
Theorem A) thanks to the time reversibility of the KdV equation.

The decay rate $\kappa $ in Theorem  \ref{main-1} has an upper
bound
\[ \kappa \leq \inf \{-\text{Re}\, \lambda : \ \lambda \in \sigma_{p }
(A_G)\} \] where $A_G$ is the operator defined by
 \[ A_G v= -v''' -\mu v - GG^* v \]
with ${\cal D} (A_G)= H^3 (\T)$ as domain.
 In order to have the decay rate $\kappa $ arbitrarily large,
 a different feedback control law is needed.

\begin{thrm}
\label{main-2}
Let $\lambda >0$, $s\geq 0$, and $\mu\in\R$  be given. There exists a
number $\delta >0$ and a linear bounded operator
 $Q_{\lambda}$ from $H^s (\T)$ to $H^s (\T)$  such that if one chooses
the feedback control law
 \[ h=-Q_{\lambda} u\]
 in system (\ref{1.1})-(\ref{1.2}), then the solution $u$ of the resulting
closed-loop system
 \begin{equation}
 \pt u+u\px u +\px ^3 u= -GQ_{\lambda} u, \quad u(x,0) =u_0 (x), \ x\in \T
\label{x-1}
 \end{equation}
 satisfies
 \[ \| u(\cdot, t) -[u_0]\|_s
\leq C e^{-\lambda t} \| u_0 -[u_0]\|_s \qquad
\text{for all } \ t\ge 0,\]
whenever $\| u_0 \|_s \leq \delta$ and $[u_0]=\mu$,
$C>0$ denoting a constant independent of $u_0$.
\end{thrm}

Note that this is  still a \emph{local} stabilization result.
However, the feedback laws in Theorems \ref{main-1} and \ref{main-2}
may be combined into a {\em time-varying} feedback law (as in
\cite{CR}) ensuring a global stabilization with an arbitrary large
decay rate.

\begin{thrm}
\label{main-3}
Let $\lambda >0$, $s\ge 0$, and $\mu\in\R$  be given. There exists a smooth map
$Q_\lambda $ from $H^s(\T )\times \R$ to $H^s(\T )$
which is periodic with respect
to the second variable, and such that the solution $u$
of the closed-loop system
$$
\partial _t u + u\partial _x u + \partial_x ^3u=-GQ_\lambda (u,t), \qquad
u(\cdot ,0)=u_0
$$
satisfies
$$
||u(\cdot , t) - [u_0]||_s \le \alpha _{s,\lambda , \mu }
(||u_0-[u_0]||_s) e^{-\lambda t}
||u_0-[u_0]||_s\qquad \text{ for all } t\ge 0,
$$
where $\alpha _{s,\lambda, \mu}: \R ^+\to \R ^+$ is a nondecreasing
continuous function depending on $s$, $\lambda$ and $\mu$.
\end{thrm}

The following remarks are in order.

\begin{rmrk} \qquad

\begin{itemize}
\item[(i)] In Theorem A, the control time $T$ is independent of the initial state $u_0$
and the terminal state $u_1$ and can be, in fact,  chosen
arbitrarily small. By contrast, in Theorem \ref{main-4}, the control
time $T$ depends on the size of the initial state $u_0$ and the
terminal state $u_1$ in the space $L^2 (\T)$. Whether the time $T$
can be chosen independent of the size of $u_0$ and $u_1$ is an
interesting open question.
\item[(ii)] While the decay rates $\kappa $ in
Theorem \ref{main-1} and $\lambda $
in Theorem \ref{main-3} are independent of $u_0$, the constants
$\alpha _{s,\mu} (\|u_0-[u_0] \|_0) $
or $\alpha_{s, \lambda ,\mu } (\| u_0-[u_0]\|_s) $ are
likely not uniformly bounded; i.e., it may happen that
\[ \lim _{r\to \infty}
\alpha _{s,\mu} (r) =\infty \ or \ \lim _{r\to \infty} \alpha _{s,
\lambda,\mu} (r) =\infty .\]
\end{itemize}
\end{rmrk}

To prove our global controllability and stabilization results
described above, we will as usual consider first the associated
linear open-loop system
\begin{equation}
\label{open} u_t +u_{xxx} = G h
\end{equation}
and the associated linear closed-loop system
\begin{equation} \label{close} u_t+u_{xxx}
= -GQ_{\lambda}.
\end{equation}
Without much difficulty we can show
by using a standard approach in control theory of linear
systems that the system (\ref{open}) is exactly controllable in the
space $H^s (\T)$ and that the closed-loop system (\ref{close}) is
exponentially stable in the space $H^s (\T)$ with an arbitrarily large
decay rate $\lambda $. However, how to extend the  linear results to
the corresponding nonlinear systems is a challenging task. Indeed,
after having published their linear results \cite{rz}, Russell and Zhang
had to wait for several years  to extend their results to
the nonlinear systems \cite{rz-1} until Bourgain
\cite{bourgain-1bis} discovered a subtle smoothing property of
solutions of the KdV equation posed on a periodic domain $\T$ when
he showed surprisingly  that the Cauchy problem of the KdV equation
(\ref{1.0}) is well-posed in the space $H^s(\T)$ for $s\geq 0$. This
then newly discovered smoothing property of the KdV equation has
played a crucial role in the proofs of Theorem A and Theorem
B in \cite{rz-1}. By contrast, establishing the global exact controllability and
stabilizability for the nonlinear system (\ref{1.4}) is even more
challenging. After all, the results presented in Theorem A and
Theorem B are essentially linear in nature; they are more or less
small perturbation of the linear results.  The global results
presented in Theorem \ref{main-4}, Theorem \ref{main-1} and Theorem
\ref{main-3} are truly nonlinear and their proofs demand new tools. The
needed help turns out to be certain propagation properties of
compactness and regularity for the KdV equation which are inspired
by those established by Laurent in \cite{laurent} for the
Schr\"odinger equation. This strategy has already been successfully applied by Dehman, Lebeau, and Zuazua \cite{DLZ} for the wave equation,
and by Dehman, G\'erard, and Lebeau \cite{DGL} and Laurent \cite{laurent,laurent2} for the Schr\"odinger equation.

Note that for any solution $u$ of the systems in consideration, its
 mean value $[u]$ is invariant.  Thus it is convenient to introduce the number
$\mu := [u]=[u_0]$, and to set
$$
\tilde u = u-\mu.
$$
Then $[\tilde u]=0$ and $\tilde u$ solves
$$
\partial _t \tilde u + \tilde u \tilde u_x + \partial _x ^3 \tilde u
 + (\mu + \tilde u ) \partial _x \tilde u = Gh.
$$
if $u$ solves (\ref{1.1}).  Throughout the paper, $\mu$ will denote
a given (real) constant, $H^s_0(\T ) = \{ u\in H^s(\T ); \ [u]=0\}$,
and $L^2_0(\T ) = \{ u\in L^2(\T ); \ [u]=0\}$. We shall establish
exponential stability results in $H^s_0(\T )$ for the equation
$$
\partial _t u + \partial _x ^3 u + \mu \partial _x u
+ u \partial _x u = - K_\lambda u
$$
that will imply all the results stated above.

\medskip
The paper is outlined as follows.

\smallskip
- In Section 2, the exact controllability and stabilizability are
presented for the associated linear systems.

\smallskip
- In Section 3,  some preliminary results in Bourgain spaces, including the
propagation of compactness and the propagation of regularity for the
KdV equation, are provided.

\smallskip
- In Section 4,  the stabilization of the KdV equation by a time invariant
feedback control law is studied.

\smallskip
- In Section 5, the stabilization of the KdV equation by a time-varying feedback
control law is investigated.

\medskip
Finally we end our introduction with a few comments on the boundary
controllability of the KdV equation posed on a finite interval
$(0,L)$:
\begin{equation} \label{boundary} \left \{
\begin{array}{l} u_t +u_x +uu_x +u_{xxx}=0, \quad x\in (0,L), \ t>0,
\\  \\u(0,t)=h_1(t),\quad
u(L,t)=h_2(t),\quad u_x(L,t)=h_3(t).\end{array}
\right.\end{equation} The problem was first investigated by Rosier
\cite{R1997} and has been intensively studied  in the past decade.
(See
\cite{R1997,Za,R2000,PMVZ,R2004,coron,Pazoto,RZ2006,LP1,GG,CC1,CC2,LP2}
and the references therein.) In contrast to control problems of
other equations (parabolic equation or hyperbolic equations for
instance), the boundary control system (\ref{boundary}) has some
interesting properties.
\begin{itemize}
\item[(i)] If
\[ L \in  {\cal N }:= \left \{ 2\pi \sqrt{\frac{j^2+l^2+jl}{3}}; j, l\in
\mathbb{N}^* \right \},\] the linear system
\[ \left \{ \begin{array}{l} u_t +u_x  +u_{xxx}=0, \quad x\in (0,L), \ t>0,
\\  \\u(0,t)=h_1(t),\quad
u(L,t)=h_2(t),\quad u_x(L,t)=h_3(t).\end{array} \right.\] associated
to (\ref{boundary}) is \textbf{not} exactly controllable if
$h_1=h_2\equiv 0$. However, the nonlinear system (\ref{boundary}) is
locally exactly controllable (still with $h_1=h_2\equiv 0$)
 \cite{R1997, coron, CC1, CC2}.
\item[(ii)] The system (\ref{boundary}) is exactly
controllable from the right (using $h_2$ or $h_3$ as control inputs
with $h_1\equiv 0$), but \textbf{only} null controllable from the
left (using $h_1$ as a control input with $h_2=h_3\equiv 0$). The
system thus behaves like a parabolic system if control is acted only
on the left end of the spatial domain  and behaves like a hyperbolic
system is control is allowed to act on the right end of the spatial
domain \cite{R2004,GG}.
\end{itemize}
\section{Linear Systems}
\setcounter{equation}{0}

 Consideration is first given to the  associate linear open loop control
system
\begin{equation}
\pt v   +\px ^3v +\mu \px v  =Gh, \quad v(x,0) = v_0 (x), \quad x\in \T,\  t\in
\R, \label{2.1}
\end{equation}
where the operator $G$ is as defined in Section 1 and $h$ is the
applied control function.

Let $A$ denote the operator
\[ Aw = -w''' -\mu w'\]
with its domain ${\cal D} (A)=H^3 (\T)$. The operator $A$ generates
a strongly continuous group $W(t)$ on the  space $L^2 (\T)$; the
eigenfunctions are simply the orthonormal Fourier basis functions in
$L^2 (\T)$,
\[ \phi _k (x) = \frac{1}{\sqrt{2\pi}} e^{ikx}, \qquad k=0, \pm 1,
\pm 2, \cdots .\] The corresponding eigenvalue of $\phi _k$ is
$$\lambda _k = ik^3 -i\mu k , \ k=0, \pm 1, \pm 2, \cdots .$$
For any $l\in \Z$, let
$$
m(l)= \# \{ k\in \Z; \ \lambda _k=\lambda _l\}.
$$
In addition,  $A^*=-A$,  $G^*=G$ and  $W^*(-t)=W(t)$ for any
$t\in \R$. Using the gap condition
$$
\lim_{|k|\to \infty} |\lambda _{k+1}-\lambda _k|=+\infty
$$
and the fact that $m(l)\le 3 $ for any $l$ and
$m(l)=1$ for $|l|$ large enough, we may deduce from Ingham lemma that
the system \eqref{2.1} is exactly controllable in $H^s_0(\T )$
in small time for any $s\ge 0$.

\begin{thrm} \cite[Theorem 2.1 and Corollary 2.1]{rz-1}
\label{linear-1}
Let $s\geq 0$ and $T>0$ be given. There exists a bounded
linear operator
\[ \Phi : H^s_0 (\T) \times H^s _0 (\T) \mapsto L^2 (0,T;
H^s_0(\T))\]
 such that for any $v_0, \ v_1 \in H^s_0 (\T)$,
\[ W(T)v_0 +\int ^T_0 W(T-t) G (\Phi (v_0, v_1)) (t)\, dt =v_1 \]
and
\[ \| \Phi (v_0, v_1)\| _{L^2 (0,T; H^s (\T))} \leq C (\| v_0\|_s
+\| v_1 \|_s ) \] where $C>0$ depends only on $T$ and $\|g\|_s$.
\end{thrm}

The following estimate  is  a direct consequence of Theorem \ref{linear-1}.
\begin{crllr} \label{est}
Let $T>0$ be given. There exists $\delta >0$   such that
\[ \int ^T_0 \| G W(t) \phi \|_0 ^2 (t) dt \geq \delta \| \phi \|_0^2 \]
for any $\phi  \in L^2 _0(\T)$. \end{crllr}

Note that the arguments presented in this paper give another proof of Corollary \ref{est}. 

In addition, if one chooses the following simple feedback
law

\[ h(v)=-G^* v,\]
the resulting closed-loop system

 \begin{equation}
 \pt v+\px^3 v +\mu \px v
=-GG^*v, \quad v(x,0) =v_0 (x), \quad x\in \T  \label{y-1}
 \end{equation}
 is exponentially stable.
 \begin{prpstn}
\label{linear-2}
  Let $s\geq 0$ be given. There exists a number
$\kappa >0$ independent of $s$ such that
 for any $v_0\in H^s_0 (\T)$, the corresponding solution
$v$ of (\ref{y-1}) satisfies
 \[ \| v (.,t)\|_s \leq C e^{-\kappa t} \| v_0 \|_s \]
 for any $t\geq 0$ where $C>0$ is  a constant depending only on $s$.
 \end{prpstn}
 \noindent
 \begin{proof} The case $s=0$ has been proved in \cite[Theorem 2]{rz}.
We only provide the proof for the
case $s=3$. The case of $0< s< 3$ follows by interpolation. The
other cases of $s$ can be proved similarly.

Pick any $v_0\in H^3_0(\T )$ and let $w=\partial _t v$. Then $w$ solves
\[ \pt w+ \px ^3w +\mu \px w
= -  G G^* w, \quad w(x,0) =w_0(x), \quad x\in \T \]
where $w_0(x)= -v^{'''}_0 (x) -\mu v_0' -  G G^* v_0 (x)$ belongs to
$L^2_0(\T )$.
Thus
\[ \| w(\cdot, t) \| _0 \| \partial _t v (\cdot, t)\| _0 \leq C_0 e^{-\kappa t} \| w_0
\|_0 \] for any $t\geq 0$. From the equation
\[ \px ^3 v + \mu \px v= -w -  G G^* v \]
it follows that
\[ \| v(\cdot,t)\|_3 \leq C_3 e^{-\kappa t} \| v_0\| _3
\]
for any $t\geq 0$. The proof is complete.
\end{proof}

\medskip
Next we show that it is possible to choose an
appropriate linear feedback law such that the
decay rate of the resulting closed-loop system is as large as one desires.

For given $\lambda > 0$, define
\[ L_{\lambda} \phi = \int ^1_0 e^{-2\lambda \tau} W(-\tau) GG^*
W^*(-\tau ) \phi \, d\tau
\]
for any $\phi \in H^s (\T)$.
Clearly, $L_{\lambda} $  is a bounded linear operator from $H^s (\T)$
to $H^s (\T )$.
Moreover, $L_{\lambda}$ is a self-adjoint positive operator on
$L^2_0(\T )$, and so is its inverse $L^{-1}_{\lambda}.$
$L_{\lambda}$ is therefore an isomorphism from $L^2_0(\T)$ onto itself.
The following result claims that the same is true on $H^s_0(\T)$.
\begin{lmm}
\label{iso}
$L_\lambda$ is an isomorphism from $H^s_0(\T)$ onto $H^s_0(\T)$ for all
$s\ge 0$.
\end{lmm}
\begin{proof}
Since the result is known for $s=0$, and $L_\lambda$ maps $H^s_0(\T)$
into itself, we only have to prove that for any $v\in L^2_0(\T)$,
$L_\lambda v \in H^s_0(\T )$ implies $v \in H^s_0(\T )$, i.e.
$D^s v\in L^2(\T )$. Using the continuity of  $L_\lambda ^{-1}$ on
$L^2_0 (\T )$ and a commutator estimate similar to
\cite[Lemma A.1]{laurent}, we obtain
\begin{eqnarray*}
||D^sv||_0
&\le& C || L_\lambda D^s v ||_0 \\
&\le& C || \int_0^1 e^{-2\lambda \tau} W(-\tau) GG^* W^* (-\tau) D^s v\,
d\tau ||_0 \\
&\le& C || D^s \int_0^1 e^{-2\lambda \tau} W(-\tau) GG^* W^* (-\tau) v\,
d\tau ||_0 \\
&&\quad +\ \
C||\int_0^1 e^{-2\lambda \tau} W(-\tau) [GG^*,D^s] W^* (-\tau)v \, d\tau||_0 \\
&\le& C||L_\lambda v||_s + C_s ||v||_{s-1}.
\end{eqnarray*}
The result follows at once for $s\in [0,1]$. An induction yields the
result for any $s\ge 0$.
\end{proof}

Choose the feedback control
\[ h = -G^*L^{-1}_{\lambda } v.\]
The resulting closed-loop system reads:
\begin{equation}
\pt v + \px ^3v  + \mu \px v
= - K_{\lambda} v , \quad v(x,0) =v_0 (x), \ x\in\T,
\label{2.2}
\end{equation}
with
\[ K_{\lambda} := GG^*L^{-1}_{\lambda }.\]
If $\lambda =0$,  we define $K_0 = GG^*$.
\begin{prpstn} \label{linear-3}
Let $s\geq 0$  and $\lambda >0$ be given. For any $v_0 \in H^s_0
(\T)$, the system (\ref{2.2}) admits a unique solution
$v\in C(\R ^+;H^s_0 (\T))$.
Moreover, there exists $M=M_s$ depending on $s$ such
that
\[ \| v(.,t)\|_s \leq M_s e^{-\lambda t} \| v_0 \| _s\] for any $t\geq 0$.
\end{prpstn}
\begin{proof}
The case $s=0$ follows from
\cite[Theorem 2.1]{slemrod}. The other cases of $s$ are proved as
for Proposition \ref{linear-2}.
\end{proof}

\section{Preliminaries}
\setcounter{equation}{0}
In this section we present some results which are essential
to establish  the exact
controllability and stabilizability of the nonlinear systems.

\subsection{The Bourgain space and its properties.}
\label{sectionprop}
For given $b,s\in \R$, and a function $u: \T \times \R \to \R$, define
the quantities
\begin{eqnarray*}
||u||_{X_{b,s}} &:=&
\left ( \sum_{k=-\infty}^{\infty}
\int_{\R}\left\langle  k\right\rangle^{2s}\left\langle
 \tau- k^3  + \mu k \right\rangle^{2b}
\left|\widehat{\widehat{u}}(k,\tau)\right|^2
 d\tau \right ) ^{\frac12}, \\
||u||_{Y_{b,s}} &:=&
\left( \sum_{k=-\infty}^{\infty}
\big( \int_{\R}\left\langle  k\right\rangle^{s}
\left\langle  \tau - k^3  + \mu k \right\rangle^{b}
\left|\widehat{\widehat{u}}(k,\tau)\right|
 d\tau \big)^2  \right ) ^{\frac12} \\
\end{eqnarray*}
where
$\widehat{\widehat{u}}(k,\tau)$
denotes the Fourier transform of $u$ with respect to the space variable $x$
and the time variable $t$ (by contrast,
$\widehat{u}(k,t)$  denotes the Fourier transform in space variable
$x$) and $\left\langle \cdot \right\rangle =\sqrt{1+|\cdot |^2}$.
Moreover, denote by $D^r$ the operator defined on $\mathcal{D'}(\Tu)$ by
\bnan
\begin{array}{rclc}
\label{Dr}
\widehat{D^r u}(k)&=&|k|^r \widehat{u}(k)
& \quad \textnormal{if}\quad k\neq 0,\\
&=& \widehat{u}(0) & \quad \textnormal{if} \quad k= 0.
\end{array}
\enan The Bourgain space $X_{b,s}$ (resp. $Y_{b,s}$)
associated to the KdV equation on $\T$
is the completion of the  space $ {\cal S}
(\T\times\R )$ under the norm $\| u\| _{X_{b,s}}$ (resp. $\| u\| _{Y_{b,s}}$).
Note that for any $u\in X_{b,s}$,
\[ \| u\|_{X_{b,s}} = \| W(-t)u\|_{H^b (\R, H^s (\T))} .\]
For given $b,s\in \R$, let
$$
Z_{b,s}=X_{b,s}\cap Y_{b-\frac{1}{2},s}
$$
be endowed with the norm
$$
||u||_{Z_{b,s}}||u||_{X_{b,s}} + ||u||_{Y_{b-\frac{1}{2},s}}.
$$
For a given  interval $I$, let
$X_{b,s} (I)$ (resp. $Z_{b,s}(I)$) be the restriction space of $X_{b,s}$ to
the interval $I$ with the norm
\begin{eqnarray*}
\| u \|_{X_{b,s} (I)}
&=& \inf \left\{ \| \widetilde{u}\| _{\Xsb} \left|\
\widetilde{u}=u \textnormal{ on } \Tu \times I \right.  \right\}\\
(\text{resp.} \ \  \| u \| _{Z_{b,s} (I)}
&=& \inf \left\{ \| \widetilde{u} \| _{Z_{b,s}} \left|\
\widetilde{u}=u \textnormal{ on } \Tu \times I  \right. \right\}).
\end{eqnarray*}
For simplicity,  we denote $X_{b,s} (I)$  (resp.
$Z_{b,s}(I)$) by $\Xsbt$  (resp. $\Zsbt$) if $I=(0,T)$.
The following properties of the spaces $X_{b,s}^T$ are $Z_{b,s}^T$
are easily verified.
\begin{itemize}
\item[(i)] $X_{b,s}
(I)$ is a Hilbert space.
\item[(ii)] $D^r u\in X_{b, s-r}(I)$ for any $u\in X_{b,s}
(I)$.
\item[(iii)]If $b_1\leq  b_2 $  and $s_1\leq s_2 $, then $X_{b_2, s_2
}$ is continuously imbedded in the space $X_{b_1, s_1}$.
\item[(iv)] For a given finite interval $I$,  if $b_1<  b_2 $  and $s_1< s_2
$, then the space $X_{b_2, s_2 } (I)$ is compactly imbedded in the
space $X_{b_1, s_1}(I)$.
\item[(v)] $Z_{\frac{1}{2},s}(I)\subset C(\overline{I};H^s(\T ))$ for any
$s\in \R$.
\end{itemize}


\begin{lmm}
\label{lin-1}
Let $b,s\in \R $ and $T>0$ be given.
There exists a constant $C>0$
such that
\begin{itemize}
\item[(i)] for any $\phi \in H^s (\T)$,
\begin{eqnarray*}
\left \| W(t) \phi \right \| _{X^T_{b,s}} &\leq& C \| \phi \| _s; \\
\left \| W(t) \phi \right \| _{Z^T_{b,s}} &\leq& C \| \phi \| _s;
\end{eqnarray*}
\item[(ii)] for any $f\in X^T_{b-1, s}$,
\[ \left \| \int  ^t_0 W(t-\tau ) f(\tau )d \tau \right \|
_{X^T_{b,s}} \leq C \| f\| _{X^T_{b-1, s}} \]
provided that $b>\frac{1}{2}$;
\item[(iii)] for any $f\in Z^T _{-\frac{1}{2},s}$,
\[
\| \int_0^t W(t-\tau ) f(\tau )\, d\tau \|
_{Z^T_{\frac{1}{2},s}} \le C\| f\|_{Z^T_{-\frac{1}{2},s}}\cdot
\]
\end{itemize}
\end{lmm}
\noindent
\begin{proof}
See e.g. \cite{tao} or \cite{CKSTT}.
\end{proof}

\begin{lmm}(Strichartz estimates)
The following estimates hold:
\begin{eqnarray}
||\sum_{k,l\in \Z}c_{k,l} e^{i(kx + lt)} ||_{L^4(\T ^2)}
&\le& C\left( \sum_{k,l\in \Z} (1+|l-k^3 + \mu k|)^{\frac{2}{3}}
|c_{k,l}|^2 \right) ^{\frac{1}{2}}, \label{E0}\\
||u||_{L^4(\T ^2 )}
&\le& C||u||_{X_{\frac{1}{3},0}}, \label{E1}\\
||u||_{L^4(\T \times (0,T))} &\le& C||u||_{X_{\frac{1}{3},0}^T}. \label{E2}
\end{eqnarray}
\end{lmm}
\begin{proof}
\eqref{E0} comes from \cite[Proposition 7.15]{bourgain-1bis}. To prove
\eqref{E1}, pick any $u\in X_{\frac{1}{3},0}$ decomposed as
$$
u(x,t)=\sum_{k\in \Z} \int_\R \widehat{\widehat u}(k,\tau)
e^{i(kx+\tau t)}d\tau .
$$
Writing $\tau = l+\sigma$ with $l\in \Z$, $\sigma \in [0,1)$, we have that
$$
u(x,t)=\int_0^1 e^{i\sigma t} \sum_{k,l\in \Z}
\widehat{\widehat u} (k,l+\sigma ) e^{i(kx+lt)}d\sigma .
$$
Using \eqref{E0} and Cauchy-Schwarz, we obtain
\begin{eqnarray*}
||u||_{L^4(\T ^2)}
&\le &
\int_0^1 || \sum_{k,l\in \Z} \widehat{\widehat u}(k,l+\sigma)
e^{i(kx+lt)}||_{L^4(\T ^2)}\, d\sigma \\
&\le & C \int_0^1
\left( \sum_{k,l\in \Z}(1+|l-k^3 +\mu k|)^{\frac{2}{3}}
|\widehat{\widehat u}(k,l+\sigma )|^2  \right)^{\frac{1}{2}}d\sigma \\
&\le& C\left( \sum_{k\in \Z}\int_0^1\sum_{l\in \Z}
(1+|l-k^3 +\mu k|)^{\frac{2}{3}}|\widehat{\widehat u}(k,l+\sigma )|^2d\sigma
\right)
^{\frac{1}{2}} \\
&\le& C
\left(\sum_{k\in \Z}\int_\R (1+ |\tau -k^3 +\mu k|)^{\frac{2}{3}}
|\widehat{\widehat u}(k,\tau )|^2 d\tau \right)^{\frac{1}{2}}.
\end{eqnarray*}
It remains to establish \eqref{E2}. Let $T>0$ and
$u\in X_{\frac{1}{3},0}^T$. Pick $p\in \N ^*$ with $T\le 2\pi p$, and
an extension $\tilde u\in X_{\frac{1}{3},0}$ of $u$ with
$||\tilde u||_{X_{\frac{1}{3},0}} \le 2 ||u||_{X_{\frac{1}{3},0}^T}$.
Then
$$
||u||^4_{L^4(\T \times (0,T))} \le ||\tilde u||^4
_{L^4(\T \times (0,2\pi p))} \le p(C ||\tilde u||_{X_{\frac{1}{3},0}})^4
\le C' ||u||^4_{X_{\frac{1}{3},0}^T}.
$$
Note that $C'$ depends only on $T$.
\end{proof}

\begin{lmm}[Bilinear estimates]
\label{lin-2} Let $s\ge 0$, $T\in (0,1)$, and
$u,v\in X^T_{\frac12 ,s}\cap L^2(0,T; L^2_0(\T ))$.
Then
there exist some constants $\theta >0$ and $C>0$ independent of $T$
and $u,v$ such that
\begin{equation}
\| (uv)_x \|_{Z^T_{-\frac12, s}}
\leq CT^{\theta } \| u\|_{X^T_{\frac12 ,s}} \| v\|_{X^T_{\frac12 ,s} }
\label{bil1}
\end{equation}
\end{lmm}

The proof of Lemma \ref{lin-2} can be found in
\cite{bourgain-1bis}  with $\theta =1/12$ (see also \cite{CKSTT}).

\bigskip

To end this section, we prove a multiplication property of the
Bourgain space $X^T_{b,s}$.  If $\psi = \psi (t)$ is any $C^{\infty}$
function, then  $\psi u \in X_{b,s}^T$
for any $u\in X^T_{b,s}$. However,
if $\phi = \phi (x)  \in C^{\infty } (\T)$, then $\phi u$ may not belong to
the space $X_{b,s}^T$ for $u\in X^T_{b,s}$. Some regularity in the
index $b$ is lost due to the fact that the multiplication by a
(smooth) function of $x$
does not keep the structure in time of the harmonics. This loss is, in fact,
unavoidable. For instance, for $k\geq 1$, let
 $u_k=\psi(t)e^{ikx}e^{i(k^3-\mu k)t}$,
where $\psi \in C^{\infty}_0(\R)$ takes the value $1$ on $[-1,1]$.
The sequence $\{ u_k\}$ is uniformly bounded
in the space $X_{b,0}$ for every $b\geq 0$.
However, multiplying $u_k$ by $\phi (x) = e^{ix}$, we observe that
$\left\|e^{ix}u_k\right\|_{X_{b,0}} \approx k^{2b}$.

The next lemma shows that this is the worst case.
\begin{lmm}
\label{lemmepseudoxsb} Let $ -1 \leq b \leq 1$, $s\in \R$ and
$\varphi \in C^{\infty}(\Tu)$. Then, for any $u\in \Xsb$,
$\varphi(x) u \in X_{b, s-2|b|}$.  Similarly, the multiplication by
$\varphi$ maps $\Xsbt$ into $X_{b,s-2|b|}^T$.
\end{lmm}
\begin{proof}
We first  consider the  case of   $b=0$  and $b= 1$.  The other cases
of $b$ will be derived later by interpolation and duality.

For $b=0$, $X_{0, s}=L^2(\R,H^s(\T))$ and the result is obvious.
For $b=1$, note that  $u\in X_{1,s}$ if and only if
 $$u \in L^2(\R,H^s(\T)) \textnormal{ and } \partial_t u+\partial_x^3 u
+\mu \partial _x u\in L^2(\R,H^s(\T)),$$
and that
 $$ \left\| u\right\|^2_{X_{1,s}}= \left\|u\right\|^2_{L^2(\R,H^s(\T))}
 +\left\|\partial_t u + \partial_x^3 u + \mu \partial _x u
\right\|^2_{L^2(\R,H^s(\T))}.$$
Thus,
\bna \left\|\varphi(x)
u\right\|^2_{X_{1,s-2}}&=&\left\|\varphi
u\right\|^2_{L^2(\R,H^{s-2}(\T))}+ \left\|\partial_t (\varphi
u)+\partial_x^3(\varphi u) + \mu \px (\varphi u)
\right\|^2_{L^2(\R,H^{s-2}(\T))}\\
&\leq& C\left(\left\|u\right\|^2_{L^2(\R,H^{s-2}(\T))}+
\left\|\varphi \left(\partial_t u+
\partial_x^3 u + \mu \px u
\right)\right\|^2_{L^2(\R,H^{s-2}(\T))}\right.\\
& &\left.+ \left\|\left[\varphi,\partial_x^3 + \mu\partial _x
\right]u\right\|^2_{L^2(\R,H^{s-2}(\T))}\right)\\
&\leq & C\left(\left\|u\right\|^2_{L^2(\R,H^{s-2}(\T))}+
\left\|\partial_t u+\partial_x^3 u + \mu \partial _x u
\right\|^2_{L^2(\R,H^{s-2}(\T))}+
\left\|u\right\|^2_{L^2(\R,H^{s}(\T))}\right)\\
&\leq & C\left\|u\right\|^2_{X_{1,s}}. \ena Here, we have used the
fact that
$$\left[\varphi,\partial_x^3 +\mu \px\right]=-3(\partial_x
\varphi)
\partial^2_x-3(\partial^2_x \varphi) \partial_x-
\partial^3_x \varphi -\mu \px \varphi $$
is a differential operator of order $2$. To conclude,
we prove that the $\Xsb$ spaces are in interpolation.
First, using Fourier transform, $X_{b,s}$ may be viewed as the
weighted $L^2$ space $L^2(\R_\tau  \times \Z _k,
\langle k\rangle ^{2s}
\langle \tau -k ^3 +\mu k \rangle ^{2b} \lambda \otimes \delta )$,
where $\lambda$ is the Lebesgue measure on $\R$ and $\delta$
is the discrete measure on $\Z$.
Then, we use the complex interpolation theorem of Stein-Weiss
for weighted $L^p$ spaces (see \cite[p. 114]{BL}):
for $0<\theta <1$
$$\left(X_{0,s},X_{1,s'}\right)_{[\theta]} \approx L^2\left(\R \times \Z,
\left\langle k\right\rangle ^{2s(1-\theta)+2s'\theta}
\left\langle \tau-k^3 + \mu k \right\rangle^{2\theta }\mu \otimes \delta
\right)\approx X_{\theta, s(1-\theta)+s'\theta}.$$
Since the multiplication by $\varphi$
maps $X_{0,s}$ into $X_{0,s}$ and $X_{1,s}$ into $X_{1,s-2}$, we
conclude that for $0\leq b \leq 1$, it maps
$\Xsb=\left(X_{0,s},X_{1,s}\right)_{[b]}$ into
\mbox{$\left(X_{0,s},X_{1,s-2}\right)_{[b]} =X_{b, s-2b}$},
which yields the $2b$ loss of regularity as announced.\\
Then, by duality, this also implies that for $0\leq b \leq 1$,
the multiplication by $\varphi(x)$ maps $X_{-b, -s+2b}$ into $X_{-b,-s}$.
As the number $s$ may take arbitrary values in $\R$,
we also have the result for $-1\leq b \leq 0$ with a loss of $-2b=2|b|$.\\
To get the same result for the restriction spaces $\Xsbt$,
we write the estimate for an extension $\tilde{u}$ of $u$,
 which yields
\bna \left\|\varphi u\right\|_{X_{b,s-2|b|}^T} \leq \left\|\varphi
\tilde{u}\right\|_{X_{b,s-2|b|}} \leq C
 \left\|\tilde{u}\right\|_{X_{b,s}}.
\ena Taking the infimum on all the $\tilde{u}$, we get the claimed
result.\end{proof}

\subsection{Propagation of compactness and regularity}
In this subsection, we present  some properties of propagation  of
compactness and regularity for  the linear differential operator $L\partial _t +\partial ^3_x + \mu \partial _x $ associated with the
KdV equation. Those propagation properties will play a key role when
studying the global stabilizability  of the KdV equation.
\begin{prpstn}
\label{thmpropagation3} Let $T>0$  and $0\leq b'\leq b \leq 1$ be
given (with $b>0$) and suppose
that $u_n\in X_{b,0}^T$  and $f_n \in X^T_{-b,-2+2b}$ satisfy
$$\partial_t u_n + \partial_x^3 u_n +\mu\partial _x u_n
=f_n $$
for $n=1, 2,...$
Assume that there exists a constant $C>0$
such that
\begin{equation}
\label{A1}
\left\| u_n \right\|_{X_{b,0}^T} \leq C \qquad \text{ for all } n\ge 1,
\end{equation}
 and that
\begin{equation}
\label{A2}
\left\|u_n\right\|_{X_{-b, -2+2b}^T}
+ \left\|f_n\right\|_{X_{-b,-2+2b}^T} +
\left\|u_n\right\|_{X_{-b',-1+2b'}^T} \rightarrow 0
\text{ as } n\to\infty .
\end{equation}

In addition, assume that for some nonempty open set $\omega \subset \T$
it holds
\[ \mbox{$u_n \rightarrow 0$ strongly in  $L^2(0,T;L^2(\omega))$.}\]
Then
\[ \mbox{$u_n \rightarrow 0$ strongly in
$L^2_{loc}((0,T);L^2(\Tu))$.}\]
\end{prpstn}
\begin{proof}
Pick $\varphi \in C^{\infty}(\Tu)$ and $\psi \in
C^{\infty}_0((0,T))$  real valued and set
\[ \mbox{ $B=\varphi(x)D^{-2}$ and
$A=\psi(t)B$}.\]
 Then $$A^*=\psi(t)D^{-2}\varphi(x).$$   For
$\varepsilon
>0$, let
$A_{\varepsilon}=Ae^{\varepsilon\partial_x^2}=\psi(t)B_{\varepsilon}$
be a regularization of $A$. Then
\begin{eqnarray*}
\alpha_{n,\varepsilon}&:=&
([A_\varepsilon , L]u_n, u_n)_{L^2(\Tu \times (0,T))}
\\  \\ &=&
([A_{\varepsilon},\partial_x^3 +\mu\partial _x]u_n,u_n)
-(\psi'(t)B_{\varepsilon}u_n,u_n).
\end{eqnarray*}
On the other hand,
\begin{eqnarray*}
\alpha_{n,\varepsilon}&=& (f_n,A^*_{\varepsilon}u_n)_{L^2(\Tu\times (0,T))}
+(A_{\varepsilon}u_n,f_n)_{L^2(\Tu\times (0,T))}
\end{eqnarray*}
since $Lu_n=f_n$ and $L^*=-L$.
By Lemma \ref{lemmepseudoxsb},
\begin{eqnarray}
\label{estimfn} \left|(f_n,A^*_{\varepsilon}u_n)_{L^2(\Tu \times
(0,T))}\right| &\leq& \|f_n\|_{X_{-b,-2+2b}^T}
\|A^*_{\varepsilon}u_n\|_{X_{b, 2-2b}^T} \nonumber\\ \nonumber
\\&\leq& \|f_n\|_{X_{-b,-2+2b}^T} \|u_n\|_{X_{b,0}^T}
\end{eqnarray}
Consequently,
\[ \lim _{n\to \infty} \sup_{0<\varepsilon\leq 1}
\left|(f_n,A^*_{\varepsilon}u_n)_{L^2( \Tu\times (0,T) )}\right|=0.\]
Similarly,  we have
\[ \lim _{n\to \infty} \sup_{0<\varepsilon\leq 1}
\left|(A_{\varepsilon}u_n,f_n)_{L^2( \Tu\times (0,T) )}\right|=0\]
and
\[ \lim _{n\to \infty} \sup_{0<\varepsilon\leq 1}
\left|(\psi'(t)B_{\varepsilon}u_n,u_n)\right|=0.\]
Thus
\[ \lim _{n\to \infty }\sup_{0<\varepsilon \leq 1}
|\alpha_{n,\varepsilon}| = 0\]
and therefore
\[ \lim _{n\to \infty} \sup_{0<\varepsilon\leq 1}
\left|([A_{\varepsilon},\partial_x^3 + \mu\partial _x]
u_n,u_n)\right|=0.\]
In particular,
  $$  \lim _{n\to \infty }
([A,\partial_x^3 + \mu \partial _x]
u_n,u_n)_{L^2(\Tu\times (0,T) )} = 0.$$
As $D^{-2}$ commutes with $\partial_x$, we have
\bna
[A,\partial_x^3 + \mu \partial _x ]
=-3\psi(t)(\partial_x \varphi) \partial^2_x D^{-2}-
3\psi(t)(\partial^2_x \varphi)\partial_xD^{-2}-\psi(t)(\partial^3_x
\varphi + \mu \partial _x \varphi ) D^{-2}.
\ena
Using the same argument as in (\ref{estimfn}),
we get
$$(\psi(t)(\partial^3_x \varphi + \mu\partial _x \varphi )D^{-2} u_n,u_n)_{
L^2(\Tu\times (0,T))} \rightarrow 0.$$
However, for the second term, the loss of
regularity is too large if we use the estimates with the same $b$.
Using the index $b'$ instead, we have
\bna (\psi(t)(\partial^2_x \varphi)\partial_xD^{-2}u_n,u_n)&\leq&
\|\psi(t)(\partial^2_x \varphi)\partial_x
D^{-2}u_n\|_{X_{b',1-2b'}^T}\|u_n\|_{X_{-b', -1+2b'}^T}  \\
&\leq& \|u_n\|_{X_{b', 0}^T}\|u_n\|_{X_{-b', -1+2b'}^T} \ena
which tends to $0$ as $n\to \infty $, by
\eqref{A1}-\eqref{A2}. Note
that $-\partial^2_x D^{-2}$ is  the orthogonal projection
on the subspace of functions with
$\widehat{u}(0)=0$. Using Rellich Theorem combined to the fact that $b>0$,
we easily see that $\widehat{u_n}(0,t)$ tends to $0$
in $L^2(0,T)$ (strongly), and hence
$$(\psi(t)(\partial_x \varphi)\widehat{u_n}(0,t),u_n)_{
L^2(\Tu\times (0,T))}
 \rightarrow 0.$$
We have thus proved that for any $\varphi\in C^{\infty}(\Tu)$ and
any $\psi \in C^{\infty}_0((0,T))$
$$(\psi(t)(\partial_x \varphi)u_n,u_n)_{L^2(\Tu\times (0,T))} \rightarrow 0.$$

Note that  a  function $\phi \in C^{\infty} (\T)$ can be written in
the form $\partial_x \varphi$  for some function $\varphi \in
C^{\infty } (\T)$  if and only if $\int_{\Tu} \phi (x) \, dx = 0$.
Thus, for any $\chi \in C^{\infty}_0(\omega)$ and any $x_0\in \Tu$, $\phi(x)\chi(x)-\chi(x-x_0)$ can be written as $\phi=\partial_x \varphi$ for
some $\varphi \in C^{\infty} (\T)$.

Since  $u_n$ is strongly convergent to $0$ in
$L^2(0,T;L^2(\omega))$,
 $$ \lim _{n\to \infty } (\psi(t)\chi u_n,u_n)_{L^2(\Tu\times (0,T))}= 0.$$
Therefore, for any $x_0\in \Tu$,
$$\lim _{n\to \infty} (\psi(t)\chi(\cdot -x_0) u_n, u_n)_{L^2(\Tu\times
(0,T))}= 0.$$
The proof is then completed by constructing a partition of
unity of $\Tu$ involving functions of the form $\chi_i(\cdot -x^i_0)$
with $\chi_i \in C^{\infty}_0(\omega )$ and $x^i_0\in \Tu$.
\end{proof}

Next we investigate the propagation of regularity  for the
operator $L= \partial _t + \partial ^3_{x} + \mu \partial _x$.
\begin{prpstn}
\label{thmpropagreg} Let $T>0$, $0\leq b<1$, $r\in \R$  and $f\in
X^T_{-b,r}$ be given.  Let   $u\in X_{b,r}^T $ be a solution of
$$\partial_t u +\partial_x^3 u + \mu \partial _x u=f.$$
If there exists a nonempty open set $\omega$ of $\T$ such that
$u\in L^2_{loc}((0,T),H^{r+\rho}(\omega))$ for some
$\rho $ with
$$0<\rho \leq \min \{1-b, \frac12\},$$
then $u\in L^2_{loc}((0,T),H^{r+\rho}(\Tu))$.
\end{prpstn}

\begin{proof} Set $s=r+\rho$ and for $ n=1,2,...$
\[ \mbox{$u_n=e^{\frac{1}{n}\partial_x^2}u=: \Xi _n u$,
 $f_n=\Xi _n f =Lu_n$}. \]
There exists a constant $C>0$
such that
$$\left\|u_n\right\|_{X_{b,r}^T}\leq C ,
\quad \left\|f_n\right\|_{X_{-b, r}^T}\leq C  \qquad \forall \, n\geq 1.
$$
Pick $\varphi \in C^{\infty}(\Tu)$ and $\psi \in C^{\infty}_0((0,T))$
as in the proof of  Proposition \ref{thmpropagation3},
 and set
\[ \mbox{$B = D^{2s-2}\varphi(x)$ and $A=\psi(t)B$.}\]
We have
\begin{eqnarray*}
& &(L u_n,A^*u_n)_{L^2(\Tu\times (0,T))}+(Au_n,Lu_n)_{L^2(\Tu\times (0,T))}\\
&=& ([A,\partial_x^3 + \mu \partial _x]u_n,u_n)_{L^2(\Tu\times (0,T))}
-(\psi'(t)Bu_n,u_n)
\end{eqnarray*}
\begin{eqnarray*}
|( Au_n,f_n)_{L^2(\Tu\times (0,T))}|
&\leq & \| Au_n\|_{X_{b,-r}^T}\|f_n\|_{X_{-b,r}^T}\\
&\leq & \| u_n\|_{X_{b,r+2\rho-2+2b}^T}\|f_n\|_{X_{-b,r}^T} \\
& \leq & C \| u_n\|_{X_{b,r}^T}\|f_n\|_{X_{-b,r}^T}\\
&\leq & C
\end{eqnarray*}
since  $r+2\rho-2+2b \leq r$.
The same estimates for the other terms
imply that
$$\left |([A,\partial_x^3 + \mu \partial _x]
u_n,u_n)_{L^2(\Tu \times (0,T) ) }\right | \leq C.$$
 Note that
$$[A,\partial_x^3+ \mu \partial _x ]=-3\psi(t)D^{2s-2}
(\partial_x \varphi)\partial^2_x-3\psi(t)D^{2s-2} (\partial^2_x
\varphi)\partial_x- \psi(t)D^{2s-2}(\partial^3_x \varphi
+\mu \partial _x \varphi )$$
and
$2s-2+1 = 2r+2\rho-1\leq 2r$. We have
\begin{eqnarray*}
&&|(\psi(t)D^{2s-2}(\partial^2_x \varphi)\partial_x
u_n,u_n)_{L^2(\Tu\times (0,T) )}| \\
&&\qquad \leq C \|\psi(t)D^{2s-2} (\partial^2_x \varphi)\partial_x
u_n\|_{L^2(0,T;H^{-r}(\T))}\|u_n\|_{L^2(0,T;H^r(\T))} \\
&&\qquad \leq C \|u_n\|_{L^2(0,T; H^r(\T))}^2 \\
&&\qquad \leq C
\end{eqnarray*}
and
\begin{eqnarray*}
&& |(\psi(t)D^{2s-2}(\partial^3_x \varphi + \mu\partial _x \varphi
)u_n,u_n)_{L^2(\Tu\times
(0,T) )}| \\
&&\qquad \leq  C \|\psi(t)D^{2s-2} (\partial^3_x \varphi
+\mu\partial _x \varphi )
u_n\|_{L^2(0,T;H^{-r}(\T))}\|u_n\|_{L^2(0,T;H^r(\T))}\\
&&\qquad \leq  C\|u_n\|_{L^2(0,T; H^r(\T))}^2 \\
&& \qquad \leq C
\end{eqnarray*}
for any $n\geq 1$. Thus
 \bnan \label{propineg1} |(\psi(t)D^{2s-2}
(\partial_x \varphi)\partial^2_x u_n,u_n)|\leq C. \enan For any
$\chi \in C^{\infty}_0(\omega)$,
\bna
&&(\psi(t)D^{2s-2} \chi ^2
\partial^2_x u_n,u_n)\\
&&\qquad  = (\psi(t)D^{s-2} \chi  \partial^2_x u_n, \chi D^s
u_n)+(\psi(t)[D^{s-2},\chi] \chi \partial^2_x u_n, D^s u_n)\\
&&\qquad = (\psi(t)D^{s-2} \chi \partial^2_x u_n, D^s \chi
u_n)+(\psi(t)D^{s-2} \chi \partial^2_x u_n,[\chi , D^s] u_n)\\
&&\qquad \quad +(\psi(t)[D^{s-2},\chi ]\chi \partial^2_x u_n, D^s u_n)
=:I_1+I_2+I_3. \ena
We infer from the assumptions that
$\chi u \in L^2_{loc}((0,T),H^s(\T))$ and that $\chi
\partial^2_x u \in L^2_{loc}((0,T),H^{s-2}(\T))$. Thus $$\chi u_n=\Xi_n
\chi u+[\chi ,\Xi_n]u$$ is uniformly bounded in
$L^2_{loc}((0,T),H^s(\T))$ by \cite[Lemma A.3]{laurent}
and the fact that $s\leq r+1$. Applying the same argument to $\chi
\partial^2_x u_n$, we obtain
\bna
|I_1|
\leq C. \ena
It follows from
\cite[Lemma A.1]{laurent} and the fact that $u\in L^2(0,T;H^r(\T))$ that
\bna
|I_2|
&\leq& C \nor{D^{r-2} \chi
 \partial^2_x u_n}{L^2(0,T;L^2(\T ))}
\nor{D^{\rho}[\chi ,D^s] u_n}{L^2(0,T;L^2(\T))}\\
&\leq&
C \nor{u_n}{L^2(0,T;H^{r}(\T))}\nor{u_n}{L^2(0,T;H^{s-1+\rho}(\T))}\leq
C. \ena
A similar bound may be obtained for $|I_3|$.
Consequently,
$$\left|(\psi(t)D^{2s-2} \chi ^2 \partial^2_x u_n,u_n)\right| \leq C$$
for any $n\geq 1.$
Then, using (\ref{propineg1}) with $\partial_x
\varphi= \chi ^2(x)-\chi ^2(x-x_0)$ yields
$$\left|(\psi(t)D^{2s-2} \chi ^2(\cdot -x_0)
\partial^2_x u_n,u_n)\right| \leq C$$
for any $n\geq 1$.
Using a partition of unity as in the proof of Proposition
\ref{thmpropagation3}, we obtain \bna
|(\psi(t)D^{2s-2}\partial^2_x u,u)| \leq C, \ena
that is
$$\int_0^T \psi (t) \big( \sum_{k\neq 0}
|k|^{2s}|\widehat{u}(k,t)|^2\big) \, dt \leq C.$$
The proof is thus complete.
\end{proof}

\medskip
\begin{crllr}
\label{corpropagreg} Let $u\in X_{\frac{1}{2},0}^T$ be a
solution of
\begin{equation}
\label{eqKdV}
\partial_t u + \partial_x^3 u + \mu \partial_x u
+ u\partial _x u= 0\textnormal{ on }
\Tu \times (0,T).
\end{equation}
Assume that $u \in C^\infty ( \omega \times (0,T))$, where
$\omega$ is a nonempty open set in $\T$.
Then $u\in C^{\infty}(\Tu \times (0,T))$.
\end{crllr}

\begin{proof}
Recall that the mean value $[u]$ is conserved. Changing $\mu$ into
$\mu+[u]$ if needed, we may assume that $[u]=0$.
We have $u\partial_x u \in X_{-\frac{1}{2},0}^T$ by
Lemma \ref{lin-2}.
It follows from Proposition \ref{thmpropagreg} that
$u\in L^2_{loc}((0,T),H^{\frac{1}{2}}(\T))$. Choose
$t_0$ such that $u(t_0) \in H^{\frac{1}{2}}(\T)$. We can then solve
\eqref{eqKdV} in $X_{\frac{1}{2},\frac{1}{2}}^T$ with the initial
data $u(t_0)$. By uniqueness of the solution in $X_{\frac{1}{2},0}^T$,
we conclude that $u\in X_{\frac{1}{2},\frac{1}{2}}^T$.
An iterated application of Proposition \ref{thmpropagreg} yields
that $u\in L^2(0,T;H^{r}(\T))$ for every $r\in \R$, and hence $u\in
C^{\infty}(\Tu \times (0,T) )$.
\end{proof}

\begin{crllr}
\label{uniciteL2}
Let $\omega$ be a nonempty open set in $\Tu$ and let
$u\in X_{\frac{1}{2},0}^T$ be a solution of
\begin{eqnarray*}
\left\lbrace
\begin{array}{rcl}
\partial_t u + \partial_x^3 u + \mu \partial _x u + u \partial_x u
&=& 0 \quad  \textnormal{ on } \Tu \times (0,T)\\
u&=&c \quad  \textnormal{ on } \omega \times (0,T)\\
\end{array}
\right.
\end{eqnarray*}
where $c\in \R$ denotes some constant. Then $u(x,t)=c$ on $\T \times (0,T)$
\end{crllr}

\begin{proof}
Using Corollary \ref{corpropagreg}, we infer that $u\in
C^{\infty}(\Tu \times (0,T))$.
It follows that  $u\equiv c$ on
$\T \times (0,T)$ by the  unique continuation property for the KdV equation
(see \cite{sautscheurer,R2004}).
\end{proof}

\section{Nonlinear systems}
\setcounter{equation}{0}
\label{sectionprop2}
In this section, we are concerned with the stability properties of
the closed loop system
\begin{equation}
\label{4.1}
\left \{ \begin{array}{l} \partial _t u  + \partial ^3 _xu +
\mu \partial _x u + u\partial _x u
= - K_{\lambda} u, \qquad x\in \T, \ 0< t < T, \\ \\ u(x,0)
=u_0(x), \qquad x\in \T ,\end{array} \right.
\end{equation}
where $\lambda \ge 0$ is a given number and $u_0\in L^2_0(\T )$.

We first check that the system is globally well-posed in the space
$H^s_0(\T)$ for any $s\geq 0$.

\begin{thrm}
\label{wellposed}
Let $\lambda \ge 0$ and $s\geq 0$ be given. Then for any $T>0$ and
any $u_0\in H^s_0 (\T)$, there exists a unique solution
$u\in Z^{T}_{\frac12 , s}\cap C([0,T];L^2_0(\T ))$ of \eqref{4.1}. Furthermore,
the following estimate holds
\begin{equation}
\label{estim1}
\|u\|_{Z^T_{\frac12,s}}  \leq \alpha_{T,s}  ( \| u_0\| _0 ) \|u_0\|_s
\end{equation}
where $\alpha _{T,s}  : \R^+ \to \R^+$ is a nondecreasing continuous
function depending only on $T$ and $s$.
\end{thrm}
\begin{proof}
We shall first establish the existence and uniqueness
of a solution $u\in Z_{\frac{1}{2},s}^T \cap L^2(0,T;L^2_0(\T ))$ of
\eqref{4.1} for $T>0$ small enough. Then we shall show
that $T$ can be taken as large as one wishes.

Let $u_0\in H^s_0 (\T)$.
Rewrite system (\ref{4.1}) in its integral form
\begin{equation}
\label{PP1}
u(t) = W(t) u_0 -\int ^t_0 W(t-\tau ) (u\px u) (\tau )d \tau -
\int ^t_0 W(t-\tau ) [K_{\lambda}u] (\tau ) d\tau
\end{equation}
where $W(t)=e^{-t(\partial _x^3 + \mu\partial _x)}$.
For given $u_0$, define the map
\[
\Gamma (v)= W(t) u_0 -\int ^t_0 W(t-\tau ) (v\px v) (\tau )d \tau -
\int ^t_0 W(t-\tau ) [K_{\lambda}v] (\tau ) d\tau .
\]
The following estimate is needed.
\begin{lmm}
\label{lem100}
For any $\varepsilon >0$ there exists a positive constant
$C(\varepsilon )$ such that
\begin{equation}
||\int_0^t W(t-\tau) [K_\lambda v](\tau ) d\tau ||_{Z_{\frac{1}{2},s}^T}
\le C( \varepsilon ) T^{1-\varepsilon} ||v||_{Z_{\frac{1}{2},s}^T}.
\label{Z1}
\end{equation}
\end{lmm}
\noindent
{\em Proof of Lemma \ref{lem100}.}
Let $v\in Z_{\frac{1}{2},s}^T$. Pick an extension of $v$ to $\T\times\R$,
still denoted by $v$, and such that
$$
||v||_{Z_{\frac{1}{2},s}} \le 2 ||v||_{Z^T_{\frac{1}{2},s}}.
$$
Pick any $\eta\in C^\infty (\R )$ with $\eta (t)=1$ for $|t|\le 1$ and
$\eta (t)=0$ for $|t|\ge 2$. By Lemma \ref{lin-1}, it is clearly sufficient
to prove that
\begin{equation}
\label{xx10}
||\eta ^2 (t/T) K_\lambda v||_{Z_{-\frac{1}{2},s}}
\le CT^{1-\varepsilon} ||v||_{Z_{\frac{1}{2},s}}.
\end{equation}
Let us first estimate $|| \eta ^2 (t/T) K_\lambda v||_{X_{-\frac{1}{2},s}}$.
We have that
$$
||\eta ^2 (t/T) K_\lambda v||_{X_{-\frac{1}{2},s}}
\le ||\eta ^2 (t/T) K_\lambda v||_{X_{\frac{-1+\varepsilon}{2},s}}
\le CT^{\frac{1-\varepsilon}{2}} ||\eta (t/T) K_\lambda v||_{X_{0,s}}
\le CT^{1-\varepsilon} ||v||_{X_{\frac{1}{2},s}}
$$
where we used \cite[Lemma 2.11]{tao} twice
and Lemma \ref{iso}. This yields also
$$
||\eta ^2 (t/T) K_\lambda v||_{Y_{-1,s}}
\le ||\eta ^2(t/T) K_\lambda v||_{X_{\frac{-1+\varepsilon}{2},s}}
\le CT^{1-\varepsilon} ||v||_{X_{\frac{1}{2},s}}.
$$
and \eqref{xx10} follows. The proof of Lemma \ref{lem100}
is complete.\qed


It follows then from Lemmas \ref{lin-1}, \ref{lin-2} and \ref{lem100} that
there exist some positive constants $\theta ,C_1,C_2$ and $C_3$ such that
\begin{eqnarray}
|| \Gamma (v)||_{Z^T_{\frac12, s}}
&\leq& C_1 \| u_0\|_s + C_2 T^\theta ||v|| ^2_{Z^T_{\frac12, s}}
+ C_3 T^{1-\varepsilon} || v||_{Z^T_{\frac12, s}} \label{Z2} \\
||\Gamma (v_1)-\Gamma (v_2)||_{Z^T_{\frac12, s}}
&\leq& C_2 T^{\theta} ||v_1+v_2||_{Z^T_{\frac12, s}}
||v_1-v_2||_{Z^T_{\frac12, s}}  +
C_3T^{1-\varepsilon} ||v_1-v_2||_{Z^T_{\frac12, s}} \label{Z3}
\end{eqnarray}
for any $v,v_1,v_2 \in Z^T_{\frac12, s}\cap L^2(0,T; L^2_0(\T ))$.
Pick $d=2C_1 \| u_0\| _s$ and $T>0$ such that
\begin{equation}
\label{Z4}
2C_2 d T^{\theta} + C_3 T^{1-\varepsilon} \le \frac{1}{2}\cdot
\end{equation}
Then
\[ ||\Gamma (v)||_{Z^T_{\frac12 , s}} \leq d \]
and
$$
\|\Gamma (v_1)-\Gamma (v_2)\|_{Z^T_{\frac12 , s}} \leq
\frac12 \| v_1 -v_2 \|_{Z^T_{\frac12 , s}}
$$
whenever
$||v  ||_{Z^T_{\frac12 , s}}\le d$,
$||v_1||_{Z^T_{\frac12 , s}}\le d$, and
$||v_2||_{Z^T_{\frac12 , s}}\le d$.
Thus the map $\Gamma$ is a contraction in the closed ball
$B_d (0)$ of $Z^T_{\frac12 , s} \cap L^2(0,T;L^2_0(\T ))$
for the $||\cdot ||_{Z^T_{\frac12 , s}}$ norm.
Its fixed point $u$ is the desired solution of (\ref{4.1}) in the space
$Z^T_{\frac12 , s}\cap L^2(0,T;L^2_0(\T ))$.
It follows from the property (v) of the Bourgain space
$Z^T_{b,s}$ recalled in the previous section that
$u\in C([0,T];H^s_0(\T ))$ with
$$
||u||_{L^\infty (0,T;H^s(\T ))}
\le C_4 || u ||_{Z^T_{\frac{1}{2},s}}
\le 2C_1C_4||u_0||_s.
$$
Let us now pass to the global existence of the solution.
Assume first that $s=0$. The solution of \eqref{4.1} satisfies
$$
||u(.,t)||^2_0 = ||u||^2_0 -
\int_0^t(GL_\lambda ^{-1}u,Gu)_0(\tau )d\tau \quad \forall t\ge 0
$$
which yields with Gronwall lemma
\begin{equation}
\label{gronwall}
||u(.,t)||^2_0 \le ||u_0||^2_0 \, e^{Ct}
\end{equation}
with $C=||G||^2||L_\lambda ^{-1}||$.
A standard continuation argument shows that \eqref{4.1} is globally
well-posed in $L^2_0(\T )$. (Note that $||u(.,t)||_0\le ||u_0||_0$ when
$\lambda =0$ and $t \geq 0$.) Next, we show that \eqref{4.1} is globally well-posed
in the space $H^3_0(\T )$. For a smooth solution $u$ of \eqref{4.1},
let $v=u_t$. Then
\begin{equation}
\label{4.10}
\left \{
\begin{array}{l}
\partial _t v  + \partial ^3 _xv +  \mu \partial _x v + \partial _x (uv)
= - K_{\lambda} v, \qquad x\in \T, \ 0< t < T, \\ \\
v(x,0) = v_0(x), \qquad x\in \T ,
\end{array}
\right.
\end{equation}
where
$$
v_0 = -K_\lambda u_0 - u_0''' - \mu u_0'-u_0 u_0'.
$$
For $T$ fulfilling \eqref{Z4}, we have
$$
||u||_{Z_{\frac{1}{2},0}^T} \le d = 2C_1 ||u_0||_0.
$$
The same computations as those leading to \eqref{Z2} yield
$$
||v||_{Z_{\frac{1}{2},0}^T}
\le C_1||v_0||_0
+(4C_1C_2 T^\theta ||u_0||_0 + C_3T^{1-\varepsilon})
||v||_{Z_{\frac{1}{2},0}^T}
$$
and hence
$$
||v||_{Z_{\frac{1}{2},0}^T} \le 2C_1||v_0||_0
$$
for $0 < T <T_1 ( || u_0 ||_0 )$, where $T_1(\cdot)$ is a continuous
nonincreasing function. Therefore,
$$
||v||_{L^\infty (0,T; L^2(\T ))} \le C_4
||v||_{Z^T_{\frac{1}{2},0}} \le C_1' ||v_0||_0
$$
for $0<T<T_1$ and $C_1'=2C_1C_4$.
>From the equation
$$
\partial _x ^3 u = -K_\lambda u - v -\mu \partial _x u - u \partial _x u,
$$
we infer that for $0<t<T<T_1$
\begin{eqnarray*}
||\partial _x ^3 u||_0
&\le&  C_7 ||u||_0 + ||v||_0
+(C_8 + ||u||_0 ) ||\partial _x u ||_{L^\infty _x}\\
&\le&
C_7 ||u||_0 + ||v||_0
+C_9(1 + ||u||_0 )
||u||_0^{\frac{1}{2}}||\partial _x^3 u||_0^{\frac{1}{2}} \\
&\le&
\frac{1}{2} ||\partial _x ^3 u||_0
+||v||_0 + C_{10} (||u||_0 + ||u||_0^3).
\end{eqnarray*}
Consequently,
$$
||u||_{L^\infty(0,T;H^3(\T ))} \le
\alpha _{T,3}(||u_0||_0)||u_0||_3
$$
for $T<T_1(||u_0||_0)$. Combined to \eqref{gronwall}, this shows
that $u\in C(\R ^+;H^3_0(\T ))$ and that \eqref{estim1} holds true
for $s=3$. A similar result can be obtained for any $s\in 3\N ^*$.
For other values of $s$, the global well-posedness follows by
nonlinear interpolation \cite{tartar, bsz-1}. The proof is complete.
\end{proof}

Next we prove a local exponential stability result when applying the
feedback law $h=-K_\lambda u$.

\begin{thrm}
\label{local}
Let $0<\lambda '<\lambda$ and $s\ge 0$ be given.
There exists $\delta >0$ such that for any $u_0\in H^s_0(\T )$ with
$||u_0||_s\le \delta $, the corresponding solution $u$ of \eqref{4.1}
satisfies
$$
|| u(.,t) ||_s \le C e^{-\lambda 't} || u_0 ||_s \qquad
\text{ for all } t\ge 0
$$
where $C>0$ is a constant independent of $u_0$.
\end{thrm}

\begin{proof}
We proceed as in \cite{RZ2007b,RZ2009}. System \eqref{4.1} can be rewritten
in an equivalent integral form
\begin{equation}
\label{integral}
u(t)=W_\lambda (t) u_0 -
\int_0^t W_\lambda (t-\tau ) (uu_x)(\tau )\, d\tau
\end{equation}
where $W_\lambda (t)=e^{-t(\partial_x ^3 +\mu\partial _x +K_{\lambda})}$.
At this point we need to extend some estimates in Lemmas
\ref{lin-1}-\ref{lin-2} for the $C^0-$group $W_\lambda (t)$.
\begin{lmm}
\label{lem10}
Let $s\ge 0$, $\lambda \ge 0$ and  $T>0$ be given. Then there exists
a constant $C>0$ such that\\
(i) for any $\phi \in H^s_0(\T )$
$$
||W_\lambda (t)\phi ||_{Z_{\frac{1}{2},s}^T} \le C ||\phi||_s.
$$
(ii) For any $u,v\in Z_{\frac{1}{2},s}^T$
$$
||\int_0^t W_\lambda (t-\tau )(uv)_x (\tau )d\tau ||_{Z_{\frac{1}{2},s}^T}
\le C||u||_{Z_{\frac{1}{2},s}^T} ||v||_{Z_{\frac{1}{2},s}^T}.
$$
\end{lmm}
\noindent
{\em Proof of Lemma \ref{lem10}:}
An application of Duhamel formula gives
\begin{equation}
\label{W1}
W_\lambda (t)\phi = W(t)\phi -\int_0^t W(t-\tau )
[K_\lambda W_\lambda (\tau )\phi ] d\tau .
\end{equation}
Using Lemma \ref{lem100}, this yields
\begin{eqnarray*}
||W_\lambda (t)\phi ||_{Z_{\frac{1}{2},s}^T}
&\le& ||W(t)\phi ||_{Z_{\frac{1}{2},s}^T} +
||\int_0^t W(t-\tau ) [K_\lambda W_\lambda (\tau ) \phi ]
d\tau ||_{Z_{\frac{1}{2},s}^T} \\
&\le& C||\phi ||_s + C T^{ 1-\varepsilon }
||W_\lambda (t)\phi||_{Z_{\frac{1}{2},s}^T}\\
\end{eqnarray*}
(i) follows at once if $T$ is small enough, say $T<T_0$.
For $T\ge T_0$, the result follows from an easy induction.
To prove (ii), we use the identity
$$
\int_0^t W_\lambda (t-\tau )f(\tau )\, d\tau
=\int_0^t W(t-\tau) f(\tau )\, d\tau
-\int_0^t W(t-\tau) K_\lambda
\left( \int_0^\tau W_{\lambda} (\tau -\sigma ) f(\sigma )\, d\sigma
\right)d\tau
$$
which gives with $f=(uv)_x$
\begin{eqnarray*}
&&||\int_0^t W_\lambda (t-\tau ) (uv)_x(\tau ) \, d\tau ||_{Z_{\frac{1}{2},s}^T} \\
&&\qquad \le ||\int_0^t W(t-\tau) (uv)_x(\tau )d\tau ||_{Z_{\frac{1}{2},s}^T}
+ ||\int_0^t W(t-\tau) K_\lambda \left(\int_0^\tau W_\lambda (\tau -\sigma )
(uv)_x(\sigma )\, d\sigma\right) d\tau  ||_{Z_{\frac{1}{2},s}^T}\\
&&\qquad \le C||u||_{Z_{\frac{1}{2},s}^T} ||v||_{Z_{\frac{1}{2},s}^T}
+CT^{1-\varepsilon}||\int_0^t W_\lambda (t-\tau) (uv)_x(\tau )\, d\tau||
_{Z_{\frac{1}{2},s}^T} \\
\end{eqnarray*}
(ii) follows again if $T$ is small enough, say $T<T_0$.
For $T\ge T_0$, the result follows from (i) and an easy induction. $\Box$\\
For given $s\ge 0$, there exists by Proposition \ref{linear-3}
some constant $C>0$ such that
$$
||W_\lambda (t) u_0||_s \le Ce^{-\lambda t}||u_0||_s \qquad
\forall t\ge 0.
$$
Pick $T>0$ such that
$$
2Ce^{-\lambda T} \le e^{-\lambda '\, T}.
$$
We seek a solution $u$ to the integral equation \eqref{integral}
as a fixed point of the map
$$
\Gamma (v) = W_\lambda (t) u_0 - \int_0^t W_\lambda (t-\tau ) (vv_x)(\tau )\,
d\tau
$$
in some closed ball $B_M(0)$ in the space $Z_{\frac{1}{2},s}^T
\cap L^2(0,T;L^2_0(\T ))$ for the $||v||_{Z_{\frac{1}{2},s}^T}$ norm.
This will be done
provided that $||u_0||_s\le \delta$ where $\delta$ is a small number
to be determined. Furthermore, to ensure the exponential stability
with the claimed decay rate, the numbers
$\delta$  and $M$ will be chosen in such a way that
$$
||u(T)||_s \le e^{-\lambda ' \, T}||u_0||_s.
$$
By Lemma \ref{lem10}, there exist some positive constants $C_1,C_2$
(independent of $\delta$ and $M$) such that
$$
||\Gamma (v)||_{Z_{\frac{1}{2},s}^T} \le C_1 ||u_0||_s +
C_2 ||v||^2_{Z_{\frac{1}{2},s}^T}
$$
and
$$
||\Gamma (v_1) -\Gamma (v_2)||_{Z_{\frac{1}{2},s}^T}
\le
C_2 ||v_1+v_2||_{Z_{\frac{1}{2},s}^T} ||v_1-v_2||_{Z_{\frac{1}{2},s}^T}.
$$
On the other hand, since $Z_{\frac{1}{2},s}^T\subset C([0,T];H^s(\T ))$, we
have for some constant $C'>0$ and all $v\in B_M(0)$
\begin{eqnarray*}
||\Gamma (v) (T)||_s &\le&
||W_\lambda (T)u_0||_s + ||\int_0^T W_\lambda (T-t)(vv_x)(\tau )d\tau||_s\\
&\le& Ce^{-\lambda T}\delta + C'M^2.
\end{eqnarray*}
Pick $\delta =C_4M^2$, where $C_4$ and $M$ are chosen so that
$$
\frac{C'}{C_4}\le Ce^{-\lambda T},\quad
(C_1C_4+C_2)M^2\le M,\quad \text{ and }\ \  2C_2M\le \frac{1}{2}.
$$
Then we have
\begin{eqnarray*}
||\Gamma (v)||_{Z_{\frac{1}{2},s}^T} &\le& M \qquad \forall v \in B_M(0),\\
||\Gamma (v_1)-\Gamma (v_2)||_{Z_{\frac{1}{2},s}^T}
&\le& \frac{1}{2} ||v_1-v_2||_{Z_{\frac{1}{2},s}^T} \qquad
\forall v_1,v_2\in B_M(0).
\end{eqnarray*}
Therefore, $\Gamma$ is a contraction in $B_M(0)$. Furthermore,
its unique fixed point $u\in B_M(0)$ fulfills
$$
||u(T)||_s =||\Gamma (u)(T)||_s\le e^{-\lambda 'T}\delta .
$$
Assume now that $0<||u_0||_s<\delta $. Changing $\delta$ into
$\delta ':=||u_0||_s$ and $M$ into $M'=(\delta '/\delta )^{\frac{1}{2}}M$,
we infer that $||u(T)||_s\le e^{-\lambda 'T}||u_0||_s$,
and an obvious induction
yields $||u(nT)||_s \le e^{-\lambda 'nT} ||u_0||_s$ for any $n\ge 0$. As
$Z_{\frac{1}{2},s}^T \cap L^2(0,T;L^2_0(\T ))
\subset C([0,T];H^s_0(\T ))$, we infer by the semigroup
property that there exists some constant $C'>0$ such that
$$
||u(t)||_s \le C'e^{-\lambda 't}||u_0||_s
$$
provided that $||u_0||_s\le \delta$. The proof is complete.
\end{proof}

The stability result presented in Theorem \ref{local} was local.
We extend it to a global stability result in the following
theorem.

\begin{thrm}
\label{global}
Assume $\lambda =0$ in (\ref{4.1}).\footnote{Recall that $K_0 = GG^*$.}
There exists a $\kappa >0$ such that for any $R_0>0$,
there exists a constant $C>0$
such that for any $u_0 \in L^2_0 (\T)$ with $$\| u_0 \|_0 \leq R_0,$$
the corresponding solution $u$ of (\ref{4.1}) (with $\lambda=0$) satisfies
\begin{equation}
\label{decay4.4}
\|u(\cdot, t)\| _0 \leq C e^{-\kappa t} \|u_0\|_0
\qquad \text{ for all } \ t  \geq 0.
\end{equation}
\end{thrm}


Theorem \ref{global} is a direct consequence of the following
observability inequality.

\begin{prpstn}
\label{propstabilisationL2}   Let  $T>0$ and  $R_0>0$ be given.
There exists a constant $\beta >1$ such that for any $u_0\in L^2_0
(\T)$ satisfying
\[ \| u_0\| _0 \leq R_0,\]
the corresponding solution $u$ of (\ref{4.1})
satisfies
\begin{equation}\|u_0 \|_0^2 \leq \beta
\int_{0}^{T} || Gu ||_0^2 (t)  dt.\label{4.2}
\end{equation}
\end{prpstn}

Indeed, if (\ref{4.2}) holds, then it follows from the energy estimate
\begin{equation}
 \|u(\cdot,t)\|_0^2 = \| u_0
\|_0^2 - \int ^t_0 ||Gu||_0^2 (\tau ) d\tau   \quad \forall
t\geq 0\label{5.2-1}
\end{equation}
that
\[ \|u(\cdot, T)\|_0^2 \leq (1-\beta ^{-1})\|u_0\|_0^2 .\]
Thus
\[ \|u(\cdot, mT)\|_0^2
\leq (1-\beta ^{-1})^m\|u_0 \|_0^2\]
which gives \eqref{decay4.4} by the semigroup property. We obtain a constant $\kappa$ independent of $R_0$ by noticing that 
for $t>c(\| u_0\| _0)$, the $L^2$ norm of $u(.,t)$ is smaller than $1$, so that we can take the $\kappa$ corresponding to $R_0=1$.\qed

\medskip

Now we present a proof of Proposition \ref{propstabilisationL2}.

\medskip
\noindent {\bf Proof of Proposition  \ref{propstabilisationL2}:} We
prove the estimate (\ref{4.2}) by contradiction. If (\ref{4.2}) is
not true, then for any $n\geq 1$, (\ref{4.1}) admits a solution
$u_n\in Z^T_{\frac12,0} \cap C([0,T];L^2_0(\T ) )$ satisfying
$$\left\|u_n(0)\right\|_0\leq R_0$$
and
\begin{equation}
\label{majorationdampedL2}
 \int_{0}^{T} ||Gu_n ||_0^2 dt < \frac{1}{n}
\left\| u_{0,n}\right\|_0^2
\end{equation}
where $u_{0,n}=u_n (0)$.
Since $\alpha_n:=\left\|u_{0,n}\right\|_0 \leq R_0$, one can choose
a subsequence of $ \{ \alpha _n\} $, still denoted by
$\{ \alpha _n \} $, such that
 $$\lim _{n\to \infty }\alpha_n = \alpha. $$
There are two possible cases: (i)  $\alpha >0$ and (ii)  $\alpha =0$.

\noindent
\underline{(i) $\alpha >0$}\\

 Note that the sequence  $\{u_n\}$ is bounded in both spaces
$L^{\infty}(0,T;L^2(\T))$ and
 $X_{\frac{1}{2},0}^T$.  By Lemma \ref{lin-2},
the sequence $\{\partial_x(u_n^2)\}$ is
bounded in the space
 $X_{-\frac{1}{2},0}^T$. On the other hand, the space
$X_{\frac{1}{2},0}^T$ is compactly imbedded in the space
$X_{0,-1}^T$. Therefore, we can extract a subsequence of
$\{ u_n \}$, still denoted by $\{ u_n\} $, such that
\begin{eqnarray*}
u_n &\to& u \quad \text{ weakly in }
X^T_{\frac{1}{2},0}, \text{  and strongly in } X_{0,-1}^T,\\
-\frac{1}{2}\partial_x(u_n^2) &\to&  f \quad \text{ weakly in }
X_{-\frac{1}{2},0}^T,
\end{eqnarray*}
where $u\in X_{\frac{1}{2},0}^T$ and $f\in X_{-\frac{1}{2},0}^T$.
Furthermore, since   $X_{\frac{1}{2},0}^T$ is continuously imbedded
in $L^4(\Tu \times (0,T))$ by \eqref{E2}, $u_n^2$ is bounded in
$L^2(\Tu \times (0,T))$. It follows that
$\partial_x(u_n^2)$ is bounded in $$L^2(0,T;H^{-1}(\T))=X_{0, -1}^T .$$
Conducting interpolation between $X_{-\frac{1}{2},0}^T$ and
$X_{0,-1}^T$,
 we obtain that $\partial_x(u_n^2)$ is bounded in
$X_{-\frac{1-\theta}{2},-\theta}^T=X_{-\frac{1}{2}
 +\frac{\theta}{2},-\theta}^T$ for $\theta \in [0,1]$. As
$X_{-\frac{1}{2}+\frac{\theta}{2},-\theta}^T
 $ is compactly imbedded in   $ X_{-\frac{1}{2},-1}^T$ for $0<\theta <1 $,
 we can extract a subsequence of $\{ u_n \}$, still denoted by
$\{ u_n\}$, such that
$-\frac{1}{2}\partial _x (u_n^2)$ converges to $f$ strongly in
$X_{-\frac{1}{2},-1}^T$.
It follows from (\ref{majorationdampedL2})   that
$$
\int_0^T ||G u_n||_0^2 \, dt   \longrightarrow
\int_0^T ||G u ||_0^2 \, dt = 0,$$
which implies that $u(x,t)=c(t)$ on $\omega \times (0,T)$
for some function $c(t)$.
Thus, passing to the limit in \eqref{4.1}, we obtain
\begin{eqnarray}
\left\lbrace
\begin{array}{rcl}
\partial_t u + \partial_x^3 u +\mu\partial _x u
&=&f \ \qquad \textnormal{ on } \T \times (0,T), \\
u&=& c(t) \qquad  \textnormal{ on } \omega \times (0,T). \\
\end{array}
\right.
\end{eqnarray}
Let  $w_n=u_n-u$ and $f_n=-\frac{1}{2}\partial_x(u_n^2)-f
-K_0 u_n$. Note first that
\begin{equation}
\label{W10}
\int_0^T ||Gw_n||^2_0 dt \int_0^T ||Gu_n||^2_0 dt +\int_0^T ||Gu||^2_0 dt
-2\int_0^T (Gu_n, Gu)_0 dt \to 0
\end{equation}
Since $w_n\to 0$ weakly in $X_{\frac{1}{2},0}^T$,  we infer from
Rellich theorem that $\int_{\T}g(y)w_n(y,t)dy\to 0$ strongly in
$L^2(0,T)$. Combined to \eqref{W10}, this yields
$$
\int_0^T\int_{\T}g(x)^2w_n(x,t)^2 dxdt \to 0.
$$
Thus
$$\partial_t w_n + \partial_x^3 w_n +\mu \partial _x w_n=f_n$$
and
$$f_n \underset{X_{-\frac{1}{2},-1}^T}{\longrightarrow}0, \quad w_n
\underset{L^2(0,T;L^2(\tilde\omega))}{\longrightarrow}0,$$
where $\tilde\omega :=\{ g > ||g||_{L^\infty (\T ) }/ 2 \}$.

Applying Proposition \ref{thmpropagation3} with $b=\frac{1}{2}$ and
$b'=0$ yields that
$$w_n \underset{L^2_{loc}((0,T);L^2(\T))}{\longrightarrow}0.$$
Consequently,   $u_n^2$ tends to $u^2$ in $L^1_{loc}((0,T);L^1(\T))$
and  $\partial_x(u_n^2)$ tends to  $\partial_x(u^2)$ in the
distributional sense. Therefore $f=-\frac{1}{2}\partial_x(u^2)$
and $u\in X_{\frac{1}{2},0}^T$  satisfies
\begin{eqnarray*}
\left\lbrace
\begin{array}{rcl}
\partial_t u + \partial_x^3 u + \mu \partial _x u +
\frac{1}{2}\partial_x(u^2)&=& 0\qquad \textnormal{ on } \Tu \times (0,T),\\
u&=& c(t)\  \quad \textnormal{ on } \omega \times (0,T).\\
\end{array}
\right.
\end{eqnarray*}
The first equation gives $c'(t)=0$ which, combined to
Corollary \ref{uniciteL2}, yields that $u(x,t)\equiv c$
for some constant $c\in \R$. Since $[u]=0$, $c=0$, and $u_n$
converges strongly to $0$ in $L^2_{loc} ((0,T), L^2 (\T))$.
We can pick some time $t_0 \in [0,T]$ such that $u_n(t_0)$ tends to $0$
strongly in $L^2(\T)$. Since
$$ \left\|u_n(0)\right\|^2_0
=\left\|u_n(t_0)\right\|^2_0 + \int_0^{t_0} ||G u_n||_0^2  dt, $$
it is inferred that $\alpha _n
=\left\|u_n(0)\right\|_0 \rightarrow 0$ which is a contradiction to
the assumption $\alpha >0$.

\medskip
\noindent
\underline{(ii) $\alpha =0$.}

Note first that $\alpha _n>0$ for all $n$.
Set $v_n=u_n/\alpha_n$ for all $n\geq 1$. Then
$$\partial_t v_n +\partial_x^3 v_n
+\mu \partial _x v_n + K_0 v_n + \frac{\alpha_n}{2} \partial_x(v_n^2)=0$$
and
\begin{eqnarray}
\label{majorenergystabL2}
\int_{0}^{T} ||G v_n||_0^2 dt < \frac{1}{n}.
\end{eqnarray}
Because of \bnan \label{vnapprox1} \left\|v_n(0)\right\|_0= 1, \enan
the sequence  $\{v_n\}$ is bounded in  both spaces
$L^{\infty}(0,T;L^2(\T))$ and $X_{\frac12, 0}^T$. Indeed, $||v_n(t)||_0$ is
a nonincreasing function of $t$, and the boundedness of
$||v_n||_{X_{\frac12,0}^T}$
for small values of $T$ follows from an estimate similar to \eqref{Z2}
(since $\alpha _n$ is bounded).
We can extract a subsequence of $\{v_n\}$, still denoted by
$\{v_n\}$, such that $v_n \to v$ weakly in the space
$X_{\frac{1}{2},0}^T$ and strongly in the spaces
$X_{-\frac{1}{2},-1}^T$ and $X_{0,-1}^T$.
Moreover, the sequence
$\{\partial _x (v_n^2)\} $ is bounded in the space $X_{-\frac12, 0}^T$,
and therefore  $\alpha_n\partial_x(v_n^2)$ tends to $0$ in the
space $X_{-\frac{1}{2},0}^T$. Finally, $\int_0^T||G v||_0^2 dt =0$.
Thus,  $v$ solves
\begin{eqnarray}
\left\lbrace
\begin{array}{rcl}
\partial_t v + \partial_x^3 v + \mu \partial_x v
&=&0 \ \ \qquad \textnormal{ on } \Tu \times (0,T) \\
v&=&c(t) \qquad  \textnormal{ on } \omega\times (0,T). \\
\end{array}
\right.
\end{eqnarray}
We infer that $v(x,t)=c(t)=c$ thanks to Holmgren Theorem, and
that $c=0$ because of $[v]=0$.

According to (\ref{majorenergystabL2})
$$\int_0^T ||G v_n||_0^2 dt\longrightarrow0$$
and so $K_{0} v_n$ converges strongly to $0$ in
$X_{-\frac{1}{2},-1}^T$. Then, an application of Proposition
\ref{thmpropagation3} as in (i) shows
that $v_n$ converges to $0$ in $L^2_{loc}((0,T),L^2(\T))$. Thus we
can pick a time $t_0\in (0,T)$ such that $v_n(t_0)$ converges
to $0$ strongly in $L^2(\T)$. Since
$$
\left\|v_n(0)\right\|^2_0
=\left\|v_n(t_0)\right\|^2_0
+ \int_0^{t_0} ||G v_n||_0^2 dt, $$
we infer from (\ref{majorenergystabL2}) that
$\left\|v_n(0)\right\|_0 \rightarrow 0$ which is a contradiction to
(\ref{vnapprox1}). The proof is complete. \qed

\medskip
Next we show that the solution $u$ of \eqref{4.1} (with $\lambda =0$)
decays exponentially
in any space $H^s(\T )$.

\begin{thrm}
\label{thm4.8}
Assume that $\lambda=0$ in (\ref{4.1}), and let $\kappa >0$ be
the infimum of the numbers $\kappa$ given respectively
in Proposition \ref{linear-2} and in Theorem \ref{global}.
Let $s\geq 0$ and let $\kappa '\in (0, \kappa )$ be given.
Then there exists a nondecreasing continuous function
$\alpha _{s, \kappa '} : \R ^+ \to \R ^+$
such that for any $u_0\in H^s_0 (\T)$, the corresponding solution
$u$ of (\ref{4.1}) satisfies
\[ \|u(\cdot, t) \|_s \leq \alpha _{s, \kappa '} (\|u_0 \|_0)
e^{-\kappa 't} \| u_0 \|_s \]
for all $t\geq 0$.
\end{thrm}

\begin{proof}
The result for $s=0$ has already been established in Theorem \ref{global} with
$\kappa '=\kappa $.
Let us consider now the case $s=3$. Pick any number $R_0>0$
and any $u_0\in H^3_0(\T )$ with
$||u_0||_0\le R_0$. Let $u$ denote the solution of
\eqref{4.1} emanating from $u_0$ at $t=0$, and let $v=u_t$.
Then $v$ solves
 \begin{equation} \pt v + \px ^3v +\mu \partial_x v + \partial_x (uv)
=-K_{0 }v, \ v(x,0)=v_0(x), \
 x\in \T, \ t>0,\label{x-100}
 \end{equation}
 where $v_0= -K_{0} u_0 -\mu u_0' - u_0 u_0' -u_0^{'''}$.
According to \eqref{estim1}  and \eqref{decay4.4}, for any $T>0$ there exists
a number $C>0$ depending only on $R_0$ and $T$ such that
\[
\|u(\cdot, t)\| _{Z^{[t, t+T]}_{\frac12,0}} \leq C e^{-\kappa t}
||u_0||_0
 \qquad \hbox{ for all } t \ \geq 0.\]
Thus, for any  $\epsilon >0$, there exists a $t^*>0$ such that
if $t\geq t^*$, one has
\[ \|u(\cdot, t)\| _{Z^{[t, t+T]}_{\frac12,0}} \leq \epsilon .\]
At this point we need an exponential stability result for the linearized
system
\begin{equation}
\pt w + \px ^3 w + \mu\px w + \px (aw) =-K_0w, \ w(x,0)=w_0(x), \
 x\in \T, \ t>0 \label{x-2}
 \end{equation}
where $a\in Z^T_{\frac12, s}\cap L^2(0,T;L^2_0(\T ))$ is a given function.
\begin{lmm}
\label{lem20}
Let $s\geq 0$ and $a\in Z^T_{\frac12, s} \cap L^2(0,T;L^2_0(\T ))
$ for all $T>0$.
Then for any $\kappa '\in (0, \kappa)$  there exist $T>0$, $\beta >0$
such that if
\[
\sup _{n\geq 1} \| a\|_{Z^{[nT, (n+1)T]}_{\frac12, s}} \leq \beta,
\]
then
\[ \| w(\cdot, t) \|_s \leq C e^{-\kappa 't} \| w_0 \| _s
\qquad \text{ for all }\ t\ge 0, \]
where $C>0$ is a constant independent of $w_0$.
\end{lmm}

\noindent
{\em Proof of Lemma \ref{lem20}:}
First, a proof similar to those of Theorem \ref{wellposed} shows that for
any $T>0$ and any $s\ge 0$, if $a\in Z^T_{\frac12,s}  \cap L^2(0,T;L^2_0(\T ))$,
 then (\ref{x-2}) admits a  unique solution $w\in Z^T_{\frac12,s}
 \cap L^2(0,T;L^2_0(\T ))$ and
\begin{equation}
\label{Z100}
\|w\|_{Z^T_{\frac12,s}} \leq \mu (\|a\|_{Z^T_{\frac12,s}})\|w_0\|_s
\end{equation}
 where $\mu: \R^+\to \R^+$ is a nondecreasing continuous function.
 Rewrite (\ref{x-2}) in its integral form
\[ w(t) = W_0 (t) w_0 -
\int ^t_0 W_0 (t-\tau )\px (aw) (\tau ) d\tau \]
where $W_0(t) = e^{ - t (\partial _x ^3 + \mu \partial _x + K_0 )}$.
Thus, for any $T>0$, by Proposition \ref{linear-2}, Lemma \ref{lem10}
and \eqref{Z100},
 \begin{eqnarray*}
  \| w(\cdot, T) \|_s & \leq  &
C_1e^{-\kappa T} \| w_0 \|_s + C_2\|
 a\|_{Z^{T}_{\frac12, s}} \|w\|_{Z^{T}_{\frac12,s}}\\ \\
 & \leq & C_1e^{-\kappa T} \| w_0 \|_s + C_2\|
 a\|_{Z^{T}_{\frac12, s}}\mu ( \|a\|_{Z^{T}_{\frac12,s}}) \| w_0
 \|_s
 \end{eqnarray*}
where $C_1 >0$ is independent of $T$ while $C_2$ may depend on $T$.
Let
\[y_n = w (\cdot, nT) \quad \text{  for } \quad n=1,2,...\]
Then, using the semigroup property of the system (\ref{x-2}),
 \[ \| y_{n+1}\|_s\leq C_1 e^{-\kappa T} \|y_n\|_s
+ C_2 \|  a \|_{Z^{[nT, (n+1)T]}_{\frac12, s}}
\mu (\|a\|_{Z^{[nT, (n+1)T]}_{\frac12,s}}) \| y_n
 \|_s\]
 for $n\geq 1$.
 Choose $T>0$ large enough  and $\beta >0$ small enough so that
 \[C_1 e^{-\kappa T}   + C_2 \beta \mu (\beta ) = e^{-\kappa 'T} \]
 Then
\[ \| y_{n+1}\|_s \leq e^{-\kappa 'T} \| y_n \|_s \]
 for any $n\geq 1$ as long as
 \[ \sup _{n\geq 1} \| a\| _{Z^{[nT, (n+1)T]}_{\frac12, s}}\leq \beta
 .\]
 Thus
 \[ \|y_{n}\|_s \le e^{-n\kappa 'T} \|y_0\|_s \]
 for any $n\geq 1$, which implies that
 \[ \|w(\cdot, t)\|_s \leq C e^{-\kappa  ' t}\| w_0\|_s \]
for all $t\ge 0$. The proof is complete. \qed

\noindent
Choose $\epsilon < \beta $, and then apply  Lemma \ref{lem20}
to (\ref{x-100}) to obtain
\[ \| v(\cdot, t) \|_0 \leq C e^{-\kappa ' (t-t^*)} \| v (\cdot, t^*) \| _0 \]
 for any $t\geq t^*$, or
 \[ \| v(\cdot, t) \|_0 \leq C_1 e^{-\kappa ' t} \| v_0 \| _0 \]
 for any $t\geq 0$, where $C_1>0$ depends only on $R_0$. It then
 follows from the equation
 \[ \px ^3u= -K_{0 }u - u\px u -\mu \px u - v \]
 and Theorem \ref{global} that
 \[ \| u(\cdot, t) \|_3 \leq C e^{-\kappa ' t} \| u_0 \| _3 \]
 for any $t\geq 0$, where $C>0$ depends only on $R_0$.

Thus the theorem has been proved for $s=0$ and $s=3$. Using the same argument for $u_1-u_2$ and $a=u_1+u_2$ for two different solutions $u_1$ and $u_2$, we obtain the Lipchitz stability estimate needed for interpolation: \[ \| (u_1-u_2)(\cdot, t) \|_0 \leq C e^{-\kappa ' t} \| (u_1-u_2)(\cdot,0) \| _0.\]
The case of $0<s<3$ follows
by interpolation. The other cases can be proved similarly.
\end{proof}

\section{Time-varying feedback law}
\setcounter{equation}{0}

In this section we prove that it is possible to design a smooth time-varying
feedback law ensuring a semiglobal stabilization with an arbitrary large
decay rate.

Let $\lambda >0$ and $s\ge 0$ be given. According to Theorem \ref{thm4.8},
there exists a number $\kappa >0$ and a nondecreasing function
$\alpha _s$ such that any solution $u$ of
\begin{equation}
\label{5.1}
\pt u + \px ^3 u + \mu \px u + u\px u =-GG^*u
\end{equation}
emanating from $u_0\in H^s_0(\T )$ at $t=0$ fulfills
\begin{equation}
\label{5.2}
||u(t)||_s \le \alpha _s (||u_0||_0)e^{-\kappa t}||u_0||_s.
\end{equation}
On the other hand, it follows from Theorem \ref{local} that for any
fixed $\lambda '\in (0, \lambda )$, any solution $u$ of
\begin{equation}
\label{5.3}
\pt u + \px ^3 u + \mu \px u + u\px u = -K_\lambda u
\end{equation}
emanating from $u_0\in H^s_0(\T )$ at $t=0$ fulfills
\begin{equation}
\label{5.4}
||u(t)||_s \le C_se^{-\lambda ' t}||u_0||_s
\end{equation}
provided that $||u_0||_s\le r_0$, for some  constant
$C_s$ and some number $r_0\in (0,1)$.
Pick any function $\theta\in C^\infty (\R ; [0,1])$
fulfilling the following properties:
\begin{eqnarray}
&&\theta (t)=1 \qquad \hbox{ for } \delta \le t \le 1-\delta \label{P1}\\
&&\theta (t)=0 \qquad \hbox{ for } 1 \le t\le 2 \label{P2} \\
&&\theta (t+2)=\theta (t) \qquad \hbox{ for all }t\in \R \label{P3}
\end{eqnarray}
where $\delta \in (0,1/10)$ is a number whose value will be specified
later. Pick a function $\rho \in C^\infty (\R ^+ ;[0,1])$ such that
\begin{equation}
\rho(r)=1 \quad \hbox{ for } r\le r_0,
\qquad \rho (r) =0 \quad \hbox{ for } r\ge 1.
\end{equation}
Let $T>0$ be given. We consider the following time-varying feedback law
\begin{eqnarray}
\label{Kut}
K(u,t) &=& \rho (||u||_s^2) \, [\theta (\frac{t}{T}) K_\lambda u
+ \theta (\frac{t-T}{T}) G G^* u]\,
+(1-\rho (||u||_s^2 ) ) G G^* u\\
&=& GG^* \, \{
\rho (||u||_s^2) \, [\theta (\frac{t}{T}) L_\lambda ^{-1} u
+ \theta (\frac{t-T}{T}) u ] + (1-\rho (||u||_s^2 ))u \}.\nonumber
\end{eqnarray}
The following result indicates that a semiglobal stabilization with
an arbitrary decay rate can be obtained.
\begin{thrm}
Let $\lambda > 0$ and let
$K=K(u,t)$ be as given in \eqref{Kut}. Pick any
$\lambda '\in (0,\lambda ) $ and any
$\lambda '' \in (\lambda ' /2, (\lambda' + \kappa )/2)$.
Then there exists a time $T_0>0$ such that for
$T>T_0$, $t_0\in \R$ and $u_0\in H^s_0(\T )$,
the unique solution of the closed-loop  system
\begin{equation}
\label{timevarying}
\pt u + \px  ^3 u + \mu \px u + u \px u = -K(u,t),\qquad u(t_0)=u_0
\end{equation}
satisfies
\begin{equation}
||u(.,t)||_s \le \gamma _s (||u_0||_s) e^{-\lambda ''(t-t_0)}||u_0||_s
\qquad \text{ for all } t\ge t_0
\end{equation}
where $\gamma _s$ is a nondecreasing continuous function.
\end{thrm}
\begin{proof}
First, proceeding as for Theorem \ref{wellposed},
we check that the system \eqref{timevarying} is globally
well-posed in $H^s_0(\T )$. Next, rough estimates for
$||u(.,t)||_s$ are
established for the times $t$ when both $K_\lambda$ and $GG^*$ are active.
\begin{lmm}
\label{lmmgronwall}
Pick any pair $(t_0,u_0)\in \R \times H^s_0(\T )$.  Then the system
\eqref{timevarying} admits a unique solution $u:\T\times [t_0,+ \infty ) \to \R$ fulfilling
$$
u\in Z_{\frac{1}{2},s}^{[t_0,t_0+T]}\cap L^2(t_0,t_0+T;L^2_0(\T)) \qquad
\text{ for all } T>0.
$$
The following a priori estimates hold true
\begin{eqnarray}
&&\text{If } \quad ||u_0||_s\le 1,\qquad
||u(.,t)||_s\le \alpha _s(1)\qquad \text{ for all } t\ge t_0;
\label{ZZ1}\\
&&\text{If } \quad ||u_0||_s > 1,\qquad
||u(.,t)||_s\le \alpha _s(||u_0||_0)
||u_0||_s\qquad \text{ for all } t\ge t_0; \label{ZZ2}\\
&&\text{If } \quad ||u_0||_s\le R,\qquad
||u(.,t)||_s\le K_s\, e^{d_s(t-t_0)} ||u_0||_s \qquad \text{ for all } t\ge t_0,
\label{ZZ3}
\end{eqnarray}
where $K_s$ and $d_s$ denote some positive constants depending only on $s$ and $R$.
\end{lmm}
\noindent
{\em Proof of Lemma \ref{lmmgronwall}:}
Let us begin with the local existence of a solution.
Pick any pair $(t_0,u_0)\in \R \times H^s_0(\T )$. It may be seen that
$$
||K(v_1,t)-K(v_2,t)||_s \le c ||v_1-v_2||_s
\qquad \text{ for all } \ v_1,v_2\in H^s_0(\T ),\ t\in\R
$$
where $c$ denotes a positive constant independent of $v_1,v_2$ and $t$. Defining the map
$$
\Gamma (v)(t) =W(t-t_0)u_0 -\int_{t_0}^t W(t-\tau)(v\partial _x v)(\tau )\, d\tau
-\int_{t_0}^t W(t-\tau) K(v(\tau ), \tau )\, d\tau ,
$$
we infer as in the proof of Theorem \ref{wellposed} that \eqref{Z2} and \eqref{Z3}
hold for all $v,v_1,v_2 \in Z_{\frac{1}{2},s}^{[t_0,t_0+\tilde{T}]}
\cap L^2(t_0,t_0+\tilde{T};L^2_0(\T))$. Moreover, the involved constants only depend on $\theta$ for its $L^{\infty}$ norm and not on $\delta$. Let $d=2C_1||u_0||_s$ and $\tilde{T}>0$ be such that
$$
2 C_2 d \tilde{T}^\theta + C_3 \tilde{T}^{1-\varepsilon} \le \frac{1}{2}\cdot
$$
Then the map $\Gamma$ is a contraction in the closed ball $B_d(0)$ of
$Z_{\frac{1}{2},s}^{[t_0,t_0+\tilde{T}]} \cap L^2(t_0,t_0+\tilde{T};L^2_0(\T))$ for
the $||v||_{Z_{\frac{1}{2},s}^{[t_0,t_0+\tilde{T}]}}$ norm. Its fixed point is the desired
solution of \eqref{timevarying}.  Note that for some constant
$C_4>0$ we have that
$$
||u||_{L^\infty (t_0,t_0+\tilde{T};H^s(\T) )}
\le C_4 ||u||_{Z_{\frac{1}{2},s}^{[t_0,t_0+\tilde{T}]}}
\le 2C_1 C_4 ||u_0||_s.
$$
Noticing that $K(u,t)= GG^*u$ for $||u||_s>1$ and
using \eqref{5.2}, we infer that the solution $u$ of \eqref{timevarying}
is defined for all $t\ge t_0$. Moreover, \eqref{5.2} yields
\eqref{ZZ1} and \eqref{ZZ2}. Let
$$d'2 C_1\max (\alpha _s(1),\alpha _s( ||u_0||_0) ||u_0||_s).$$
Note that $d'$ depends only on $R$ and $s$.
Replacing $\tilde{T}$ by $T'$ satisfying
$$
2C_2 d'{T'}^\theta +C_3 {T'}^{1-\varepsilon} \le \frac{1}{2}
$$
in the application of the contraction mapping principle,
we infer that
the (unique) solution $u$ of \eqref{timevarying} fulfills
$$
||u||_{Z_{\frac{1}{2},s}^{[t_0+kT',t_0+(k+1)T']}}
\le 2C_1 ||u(.,t_0+kT')||_s.
$$
This gives
$$
||u||_{L^\infty (t_0+kT',t_0+(k+1)T',H^s(\T))}
\le 2C_1 C_4 ||u(.,t_0+kT')||_s
$$
and
$$
||u(.,t)||_s\le  K_s e^{d_s(t-t_0)} ||u_0||_s
$$
for some constants $K_s>0$, $d_s>0$ depending only on $s$ and $R$. \qed

Given $\lambda ''$ as in the statement of the theorem, we pick $\delta >0$ such that
\begin{equation}
\label{AB1}
\lambda '' < - 2\delta d_s + (1-2\delta )
\frac{\kappa +  \lambda '}{2} \qquad \text{ and }\qquad
\delta d_s - (1 - 2\delta ) \kappa < 0 \cdot
\end{equation}
Next, choose $r_1\in (0,r_0)$ such that
\begin{equation}
\label{AB2}
\alpha _s (\alpha _s (1))C_sK_s^4e^{4\delta Td_s} r_1 <r_0,
\end{equation}
and $T_0>0$ such that
\begin{eqnarray}
\label{AB3}
\alpha _s (1) \alpha _s (\alpha _s (1))K_s
 e^{ [\delta d_s - (1-2\delta ) \kappa ] T} &\le &r_1, \\
\label{AB4}
\alpha _s (1)C_sK_s^4e^{[4\delta d_s - (1-2\delta) (\kappa + \lambda ') ]T}
&\le&  e^{-2\lambda ''T}
\end{eqnarray}
for all $T\ge T_0$. Note that $T_0$ exists by \eqref{AB1}.
Pick any $u_0\in H^s_0(\T )$ and any time
$t_0\in \R$.
The proof rests on a series of claims. \\
{\sc Claim 1}. There exists a time
$t_1\in [t_0,t_0+\kappa ^{-1}\ln (\alpha _s (||u_0||_0) ||u_0||_s)]$
such that
\begin{equation}
\label{RR1}
||u(t_1)||_s \le 1.
\end{equation}

Without loss of generality we may assume that
$||u(t_0)||_s\ge 1$. Then the dynamics of $u$ is governed by
\eqref{5.1} as long as $||u(t)||_s\ge 1$. By \eqref{5.2}, we have
\eqref{RR1} for some time $t_1$ with
$$
\alpha _s(||u_0||_0) e^{-\kappa (t_1-t_0)}||u_0||_s \le 1.
$$
Therefore, Claim 1 holds.\qed

\noindent
{\sc Claim 2.} There exists a time $t_2\in 2\Z T\cap [t_1,t_1+3T]$
such that
\begin{equation}
\label{RR2}
||u(t_2)||_s\le r_1.
\end{equation}

\, From the fact that $||u(t_1)||_0 \le 1$ and \eqref{5.2}
we have that
$$
||u(t)||_s\le \alpha _s (1) \qquad \hbox{ for all } t\ge t_1.
$$
Pick $R=\alpha _s(1)$ and let $K_s$ and $d_s$ be as given in Lemma \ref{lmmgronwall} for that choice of $R$.
Let $t_1'\ge t_1$ denote the first time of the form
$t_1'=(2k+1)T + \delta$ with $k\in \Z$, and let $t_2=(2k+2)T$.
Then it follows from \eqref{5.2}, \eqref{ZZ3} and \eqref{AB3}
that
$$
||u(t_2)||_s\le K_s e^{\delta T d_s} \alpha _s(\alpha _s(1))
e^{-\kappa (1-2\delta )T}||u(t_1')||_s \le r_1.
$$
\noindent
{\sc Claim 3.} $||u(t)||_s\le r_0$ for all $t\ge t_2$ and
$||u(t_2+2kT)||_s \le e^{-2k\lambda ''T} ||u(t_2)||_s$
for all $k\in \N$.

First, we notice that the dynamics of $u$ is governed by \eqref{5.3}
(resp. by \eqref{5.1})  when $t\in (t_2+\delta T, t_2+(1-\delta)T)$
(resp. when $t\in (t_2+(1+\delta) T, t_2+(2-\delta )T )$), as long as
$||u (t)||_s\le r_0$. Therefore, using \eqref{5.2}, \eqref{5.4},
\eqref{ZZ3}, and \eqref{AB2} we obtain that
$$
||u(t)||_s\le ( \alpha _s (\alpha _s ( 1 ))
K_s^2e^{2\delta Td_s})(C_sK_s^2e^{2\delta Td_s})
||u(t_2)||_s <r_0\qquad \text{ for all }t\in [t_2,t_2+2T].
$$
On the other hand, by \eqref{AB4},
\begin{eqnarray*}
||u(t_2+2T)||_s
&\le& (\alpha _s(1)e^{-\kappa (1-2\delta )T}K_s^2 e^{2\delta T d_s})
(C_s e^{-\lambda ' (1-2\delta )T}K_s^2e^{2\delta T d_s})||u(t_2)||_s \\
&\le& e^{-2\lambda ''T} ||u(t_2)||_s \\
&\le& r_1.
\end{eqnarray*}
The claim follows by an obvious induction. \qed

It follows from Claim 3 that for $t\ge t_2$
$$
||u(t)||_s \le c \, e^{-\lambda ''(t-t_2)} ||u(t_2)||_s
$$
for some constant $c$ independent of $t$ and $u_0$. Since
$$
t_2-t_0\le 3T+ \kappa ^{-1}\ln (\alpha _s (||u_0||_0)||u_0||_s),
$$
the theorem follows.
\end{proof}

\begin{rmrk}

\quad

\begin{itemize}
\item A natural idea to combine both feedback controls
would be to consider a discontinuous
feedback control which agrees with $K_0 u$ when $||u||_s$ is large, and
with $K_\lambda u$ when $||u||_s$ is small. The main difficulty is
then to define properly what we mean by a solution of the
closed-loop system. In finite dimension, the Filippov solutions
are widely used by the control community to deal with discontinuous systems.
(See \cite{CR} for the definition of a Filippov solution.)
The main advantage of the time-varying feedback law considered here is
its regularity, which guarantees the existence and uniqueness of
``classical'' solutions for the closed-loop system.
\item It would be interesting to see whether a smooth time-invariant
feedback law ensuring a semi-global exponential stabilization with an
arbitrary decay rate can be designed.
\item A simpler, but less efficient, time-varying feedback law
is
$$
K(u,t) := \theta (\frac{t}{T})
\rho (||u||_s^2) \, K_\lambda u
+ \theta (\frac{t-T}{T}) G G^* u.
$$
\end{itemize}
\end{rmrk}

\bigskip
\textbf{Acknowledgments}:  This  work is partially conducted while
the third author (BZ) was visiting  \emph{Laboratoire Jacques-Louis
Lions} of Universit\'e Pierre et Marie Curie in  May and June of
2008, and  the \emph{Yantz Center of Mathematics} of Sichuan
University in November of 2008. He thanks both institutions for
their hospitality. BZ was partially supported by the \emph{Paris
Foundation of Mathematics}, the \emph{Charles Phelps Taft Memorial
Fund}  at the University of Cincinnati and the \emph{Chunhui
program} (State Education Ministry of China) under grant
Z007-1-61006


\addcontentsline{toc}{section}{References}
\begin{thebibliography}{9}

\bibitem{BL}  Bergh, J. and  L\"ofstrom, J., \emph{Interpolation Spaces,
An Introduction,} Springer Verlag 1976


\bibitem{bsz-1} Bona, J. L., Sun, S. M. and Zhang, B.-Y.,
The initial-boundary value problem for the Korteweg-de Vries
equation in a quarter plane, \emph{Trans.  American Math. Soc.} {\bf
354}(2001), 427--490.


\bibitem{bourgain-1bis} Bourgain, J.,  Fourier transform restriction
phenomena for certain lattice subsets and applications to non-linear
evolution equations, part II: the KdV  equation, {\em Geom. \&
Funct. Anal. } {\bf 3}\,(1993), 209 -- 262.

\bibitem{bouss} Boussinesq, J.,  \emph{Essai sur la th\'eorie des eaux
courantes}; M\'emoires pr\'esent\'es
par divers savants \`a l'Acad. des Sci. Inst. Nat. France,
\textbf{23} (1877), 1C680.

\bibitem{CC1} Cerpa, E.,  and Cr\'epeau, E., {\em Boundary controllability
for the nonlinear Korteweg-de Vries equation on any critical domain}, Ann.
I.H. Poincar\'e - AN (2008), doi:10.1016/j.anihpc.2007.11.003.

\bibitem{CC2} Cerpa, E., and Cr\'epeau, E., {\em Rapid exponential
stabilization for a linear Korteweg-de Vries equation}, Discrete Contin.
Dyn. Syst. Ser. B, {\bf 11}\,(2009), no. 3, 655--668.

\bibitem{CKSTT} Colliander J., Keel M., Staffilani G.,
Takaoka H., and Tao T., {\em Sharp global well-posedness for KdV and
modified KdV on $\R$ and $\T$}, Journal of the AMS {\bf 16} (2003), No. 3,
705--749.

\bibitem{coron} Coron, J.-M. and Cr\'{e}peau, E., {\em Exact boundary
controllability of  a nonlinear KdV equation with a critical
length}, J. Eur. Math. Soc., {\bf 6}\,(2004), 367--398.

\bibitem{CR} Coron, J.-M. and Rosier, L., {\em
A Relation Between Continuous Time-Varying and Discontinuous
Feedback Stabilization}, Journal of Mathematical Systems,
Estimation, and Control, {\bf 4}\,(1994), 67--84.

\bibitem{DGL} Dehman, B., G\'erard, P., Lebeau, G.,
Stabilization and control for the nonlinear Schr\"odinger equation
on a compact surface , \emph{Math. Z} {\bf 254}\,(2006), 729--749.
\bibitem{DLZ} B. Dehman, G. Lebeau and E. Zuazua : Stabilization and control for the subcritical semilinear wave equation. Anna. Sci. Ec. Norm. Super. 36:525-551 (2003).
\bibitem{soliton}  Gardner, C. S., Greene, J. M.,  Kruskal, M. D. and
Miura, R. M., Method for solving the Korteweg-de Vries equation,
\emph{Phys. Rev. Lett.,} \textbf{19} (1967), 1095--1097.

\bibitem{GG} Glass, O., and Guerrero, S., Some exact
controllability results for the linear KdV equation and uniform
controllability in the zero-dispersion limit, \emph{Asymptot. Anal.}
{\bf 60} (2008), no. 1-2, 61--100.

\bibitem{jager} de Jager, E. M., On the origin of the Korteweg-de Vries
equation, \emph{arXiv:math.HO/0602661}

\bibitem{KappelerTopalov} Kappeler, T. and Topalov, P., {\em Well-posedness
of KdV on $H^{-1}(\T )$}, Duke Math. J., {\bf 135} (2006), 327--360.

\bibitem{kato} Kato, T., On the Cauchy problem for the (generalized)
Korteweg-de Vries equations, in \emph{Advances in Mathematics Supplementary
Studies, Stud. Appl. Math.} \textbf{8,} Academic Press, New York,
1983,  93--128.

\bibitem{kpv-1} Kenig, C. E.,  Ponce, G.  and Vega, L.,  Well-posedness of the
initial value problem for the KdV equation, \emph{J. Amer. Math.
Soc.,} \textbf{4} (1991),  323--347.


\bibitem{kpv-2}  Kenig, C. E.,  Ponce, G.  and Vega, L., A
bilinear estimate with applications to the KdV equation, \emph{J.
Amer. Math. Soc.}, {\bf 9} (1996), 573--603.

\bibitem{kdv}  Korteweg, D. J. and  de Vries, G.,  On the change
of form of long waves advancing in a
rectangular canal, and on a new type of long stationary waves,
\emph{Philos. Mag.}, \textbf{39} (1895),  422--443.

\bibitem{laurent} Laurent, C.,
Global controllability and stabilization for the nonlinear
Schr\"odinger equation on an interval, \emph{ESAIM
Control Optim. Cal. Var.}, in press.

\bibitem{laurent2} Laurent, C.,
Global controllability and stabilization for the nonlinear
Schr\"odinger equation on some compact manifolds of dimension 3,
submitted.

\bibitem{LP1} Linares, F., and Pazoto, A. F.,   On the exponential decay
of the critical generalized Korteweg-de Vries equation with
localized damping, \emph{Proc. Amer. Math. Soc.} {\bf 135} (2007),
no. 5, 1515--1522.

\bibitem{LP2} Linares, F., and Pazoto, A. F.,
Asymptotic behavior of the Korteweg-de Vries equation posed in a
quarter plane, \emph{J. Diff. Equations} {\bf 246} (2009) 13421353.


\bibitem{MORZ} Micu, S., Ortega, J., Rosier, L., and Zhang, B.-Y.,
  Control and stabilization of a family of Boussinesq systems,
\emph{Discrete and Continuous Dynamical Systems,} {\bf 24} (2009),
no. 2, 273--313.

\bibitem{miura}   Miura,  R. M., The Korteweg-de Vries equation: A survey of
results, \emph{SIAM Rev.,} \textbf{18} (1976),   412--459.


\bibitem{Pazoto} Pazoto, A. F.,   Unique continuation and decay for the
Korteweg-de Vries equation with localized damping, \emph{ESAIM
Control Optim. Calc. Var., }{\bf 11} (2005), pp. 473--486.


\bibitem{PR} Pazoto, A. F. and Rosier, L., Stabilization of a
Boussinesq system of KdV-KdV type, \emph{Systems $\&$ Control
Letters} {\bf 57} (2008), 595--601.


\bibitem{PMVZ} Perla-Menzala, G., Vasconcellos, C. F.,  and Zuazua, E.,
 Stabilization of the Korteweg-de Vries equation with localized
damping, \emph{Quart. Appl. Math.,} {\bf 60} (2002), pp. 111--129.


\bibitem{ru-1} Russell, D. L, D. L.  Computational study of
the Korteweg-de Vries
 equation with localized control
action, \emph{Distributed Parameter Control Systems: New Trends and
Applications}, G. Chen, E. B. Lee, W. Littman, and L. Markus, eds.,
Lecture Notes in Pure and Appl. Math., vol. 128, Marcel Dekker, New
York, 1991, 195--203.

\bibitem{R1997} Rosier, L.,  Exact boundary controllability for the
Korteweg-de Vries equation on a bounded domain, \emph{ESAIM
Control Optim. Cal. Var.,}
{\bf 2} (1997), 33--55.

\bibitem{R2000} Rosier, L., Exact boundary controllability for the linear
Korteweg-de Vries equation on the half-line, \emph{SIAM J. Control
Optim.} {\bf 39} (2000) 331--351.

\bibitem{R2004} Rosier, L., Control of the surface of a fluid
by a wavemaker, \emph{ESAIM Control Optim. Cal. Var.}
{\bf 10} (2004), 346--380.

\bibitem{RZ2006}  Rosier, L. and  Zhang, B.-Y.,  Global
stabilization of the generalized Korteweg-de Vries equation,
\emph{SIAM J. Control Optim. } {\bf 45}\,(2006), no. 3, 927--956.

\bibitem{RZ2007b}  Rosier, L. and  Zhang, B.-Y.,  Exact
controllability and stabilization of the nonlinear  Schr\"odinger
equation on a bounded interval, \emph{SIAM J. Control Optim.} {\bf
48}\, (2009), no. 2, 972--992.

\bibitem{RZ2008} Rosier, L. and Zhang, B.-Y.,  Exact boundary
controllability of the nonlinear Schr\"odinger equation, \emph{J.
Differential Equations} {\bf 246} (2009), 4129--4153.

\bibitem{RZ2009} Rosier, L. and Zhang, B.-Y., Control and stabilization
of the nonlinear Schr\"odinger equation on rectangles, submitted.

\bibitem{rz} Russell, D. L. and Zhang, B.-Y.,  Controllability and
stabilizability of the third order linear dispersion equation on a
periodic domain,\emph{ SIAM J. Cont. Optim.,} {\bf 31}\,(1993),
659--676.

\bibitem{rz-1} Russell, D. L. and Zhang, B.-Y., Exact
controllability and stabilizability of the Korteweg-de Vries
equation, \emph{Trans. Amer. Math. Soc.,}  {\bf 348 }\,(1996),
3643--3672.

\bibitem{st}  Saut, J.-C. and Temam, R.,  Remarks on the Korteweg-de Vries
equation, \emph{Israel J. Math.,} \textbf{24} (1976), 78--87.


\bibitem{sautscheurer} Saut, J.-C. and
Scheurer, B.,   Unique continuation for some evolution equations,
\emph{J. Differential Equations}, {\bf 66} (1987), no. 1, 118--139.

\bibitem{slemrod} Slemrod, M.,
  A note on complete controllability and stabilizability
for linear control systems in Hilbert space, \emph{ SIAM Control,}
{\bf 12}\,(1974), 500-508.

\bibitem{tao} Tao, T., {\em Nonlinear dispersive equations, Local and global
analysis}. CBMS Regional Conference Series in Mathematics, 106.
AMS, Providence, RI, 2006.

\bibitem{tartar} Tartar, L.,
 Interpolation non lin\'eaire et r\'egularit\'e, J. Funct. Analys,
 {\bf 9}(1972), 469-489.

\bibitem{temam} Temam, R., Sur un probl\`eme non
lin\'eaire, \emph{J. Math. Pures Appl}.,
\textbf{48}(1969), 159--172.

\bibitem{Za-2}  Zhang, B.-Y., Unique continuation for the Korteweg-De Vries
equation, \emph{SIAM J. Math. Anal.,} \textbf{23} (1992), 55--71.

\bibitem{Za-3} Zhang, B.-Y.,  Taylor series expansion for solutions of the
Korteweg-de Vries equation with respect to their initial values,
\emph{J. Funct. Anal.,} \textbf{129 }(1995), 293--324.

\bibitem{Za} Zhang, B.-Y.,  Exact boundary controllability of the
Korteweg-de Vries equation, \emph{SIAM J. Cont. Optim}., 37 (1999),
543--565.


\end{thebibliography}
\end{document}